\documentclass{nyj}

\usepackage{graphicx}
\usepackage{amssymb}
\usepackage{verbatim}

\title[Spectral Convergence]{Spectral Convergence of the Discrete Laplacian on Models of a Metrized
Graph}

\author{X.W.C. Faber}

\address{Department of Mathematics\\
Columbia University\\
New York, NY  10027\\
USA}
\email{xander@math.columbia.edu}



\newtheorem{thm}{Theorem}
\newtheorem{lem}{Lemma}
\newtheorem{cor}{Corollary}

\newtheorem{prop}{Proposition}

\newcommand{\Dir}{\operatorname{Dir}} 
\newcommand{\Funct}{\operatorname{Funct}} 
\newcommand{\Ker}{\operatorname{Ker}} 
\newcommand{\Zh}{\operatorname{Zh}} 
\newcommand{\dist}{\operatorname{dist}} 
\newcommand{\CPA}{\operatorname{CPA}} 
\newcommand{\vect}{\operatorname{Vec}} 

\newcommand{\ip}[2]{\left< #1,#2 \right>} 
\newcommand{\CC}{\mathbb{C}} 
\newcommand{\NN}{\mathbb{N}}  
\newcommand{\p}{\varphi_N} 
\newcommand{\HH}{\mathcal{H}} 
\newcommand{\TT}{\mathcal{T}} 
\newcommand{\SSS}{\mathcal{S}} 

\begin{document}

\begin{abstract}
A metrized graph is a compact singular $1$-manifold endowed with a metric. A given
metrized graph can be modelled by a family of weighted combinatorial graphs. If one
chooses a sequence of models from this family such that the vertices become uniformly
distributed on the metrized graph, then the $i$th largest eigenvalue of the Laplacian
matrices of these combinatorial graphs converges to the $i$th largest eigenvalue of the
continuous Laplacian operator on the metrized graph upon suitable scaling. The
eigenvectors of these matrices can be viewed as functions on the metrized graph by linear
interpolation. These interpolated functions form a normal family, any convergent
subsequence of which limits to an eigenfunction of the continuous Laplacian operator on
the metrized graph.
\end{abstract}
\maketitle \tableofcontents

\section{Introduction}
Roughly speaking, a metrized graph $\Gamma$ is a compact singular $1$-manifold endowed
with a metric. A given metrized graph can be modelled by a family of weighted
combinatorial graphs by marking a finite number of points on the metrized graph and
declaring them to be vertices. On each of these models we have a Laplacian matrix, which
acts on functions defined on the vertices of the model. On the metrized graph there is a
measure-valued Laplacian operator. The goal of this paper is to prove that the spectra of
these two operators are intimately related. Indeed, we will show that if we pick a
sequence of models for $\Gamma$ whose vertices become equidistributed, then the
eigenvalues of the discrete Laplacians on the models converge to the eigenvalues of the
Laplacian operator on $\Gamma$ provided we scale them correctly. Moreover, we can show
that in a precise sense the eigenvectors of the discrete Laplacian converge uniformly to
eigenfunctions of the Laplacian operator on the metrized graph.

This type of approximation result for Laplacian eigenvalues dates back at least to the
papers of Fukaya (\cite{KeFu}) and Fujiwara (\cite{KoFu1}, \cite{KoFu2}). In \cite{KeFu}
the measured Hausdorff topology is defined on closed Riemannian manifolds of a fixed
dimension subject to certain curvature hypotheses, and it is shown that convergence of
manifolds in this topology implies convergence of eigenvalues for the associated
Laplace-Beltrami operators. An analogue of the measured Hausdorff topology for the class
of finite weighted graphs is defined in \cite{KoFu1}; a similar convergence result for
eigenvalues of the discrete Laplacian operator is obtained. In \cite{KoFu2}, Fujiwara
approximates a closed Riemannian manifold by embedded finite graphs and proves that the
eigenvalues of the Laplace-Beltrami operator can be bounded in terms of the eigenvalues
of the associated graph Laplacians. In the present paper we use an approach remarkably
similar to the one in \cite{KoFu1}; this is purely a coincidence as the author was
unaware of Fujiwara's work until after giving a proof of the main theorem. It would be
interesting to see if the methods employed here can improve upon \cite{KoFu2} when
applied to the approximation of Riemannian manifolds by graphs.

\begin{figure}[!ht]
  \scalebox{1}{\includegraphics{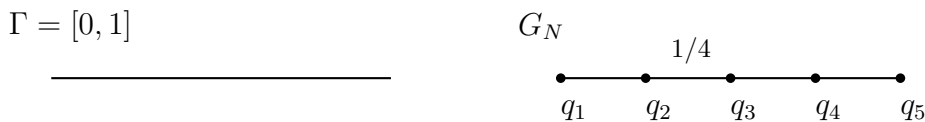}}

  \caption{Here we see the metrized graph $\Gamma$ and the graph $G_N$ where $N=5$.
           Each edge of $G_N$ has length $1/4$ or weight $4$. One can view $G_N$ as a
           discrete approximation to the metric space $\Gamma$.}

  \label{Interval}
\end{figure}

Without giving too many of the details, let us briefly display some of the content of the
main theorem (Theorem~\ref{Main Theorem}) in a special case. See section~\ref{Notation
Section} for precise definitions of all of the objects mentioned here. Let $\Gamma =
[0,1]$ be the interval of length~1. Let $G_N$ be a weighted graph with vertex set $V_N =
\{q_1, \ldots, q_N\}$, edge set $E_N = \{q_iq_{i+1}: i=1, \ldots, N-1\}$, and weight
$N-1$ on each edge. The reciprocal of the weight of an edge will be its length. We say
that $G_N$ is a model of our metrized graph $\Gamma$; see Figure~\ref{Interval}.

For each $N$ we can define the Laplacian matrix $Q_N$ associated to the graph $G_N$. It
is an $N \times N$ matrix which contains the weight and incidence data for the graph. Let
$\lambda_{1,N}$ be the smallest positive eigenvalue of $Q_N$. The following table gives
the value of $N\lambda_{1,N}$ (the scaled eigenvalue) for several choices of $N$. In each
case, the eigenspace associated to this eigenvalue has dimension~1.

\begin{center}
\begin{tabular}{|c|c|}
\hline
$N$  & $N\lambda_{1,N}$ (approximate) \\
\hline
5 &  7.6393\\
\hline
10 &  8.8098\\
\hline
50 &  9.6690  \\
\hline
100 &  9.7701 \\
 \hline
200 &  9.8201 \\
 \hline
500 &  9.8498 \\
 \hline
\end{tabular}
\end{center}

Now let $f: \Gamma \to \CC$ be a continuous function that is smooth away from a finite
set of points and that has one-sided derivatives at all of its singular points. The
Laplacian of such a function $f$ is defined to be
\[ \Delta f = -f''(x)dx - \sum_{\substack{\text{$p$ a singular} \\ \text{point of $f$}}} \sigma_p(f)
\ \delta_p,\] where $dx$ is the Lebesgue measure on the interval, $\delta_p$ is the point
mass at $p$, and $\sigma_p(f)$ is the sum of the one-sided derivatives of $f$ at the
singular point $p$. We say that a nonzero function $f$ is an eigenfunction for $\Delta$
with respect to the measure $dx$ if the following two conditions obtain:
\begin{itemize}
\item $\int_{\Gamma} f(x)dx=0$
\item There exists an eigenvalue $\lambda >0$ such that $\Delta f = \lambda f(x) dx$.
\end{itemize}
The eigenvalues of $\Delta$ with respect to the measure $dx$ are $n^2\pi^2$ for $n = 1,
2, 3, \ldots$, and the eigenspace associated to each eigenvalue has dimension~1. Denoting
the smallest eigenvalue by $\lambda_1(\Gamma)$, we see that $\lambda_1(\Gamma) =\pi^2
\approx 9.8696$. The values of $N\lambda_{1,N}$ in the above table could reasonably be
converging to $\lambda_1(\Gamma)$.

Part (A) of Theorem~\ref{Main Theorem} asserts that the scaled eigenvalues
$N\lambda_{1,N}$ do indeed converge to $\lambda_1(\Gamma)$. A similar statement can be
made for the second smallest eigenvalues of $Q_N$ and $\Delta$, as well as the third,
etc. The fact that the dimensions of the eigenspaces for $\lambda_{1,N}$ and
$\lambda_1(\Gamma)$ agree is no coincidence, and a precise statement of this phenomenon
constitutes part (B) of Theorem~\ref{Main Theorem}. If we choose an $\ell^2$-normalized
eigenvector $h_N$ of the matrix $Q_N$ for each $N$, then the sequence $\{h_N\}$ can be
viewed as a family of piecewise affine functions on the interval by linear interpolation.
Part (C) of the main result states that this family is normal, and any subsequential
limit of it will be an $L^2$-normalized eigenfunction for $\Delta$ with respect to the
$dx$ measure.

We make all of this precise in the next section after defining some of the necessary
notation and conventions. We will follow quite closely the treatment of metrized graphs
given in \cite{BR}. Applications of metrized graphs to other areas of mathematics and the
physical sciences can be found in \cite{Ku} and \cite{BR}. For a more conversational
approach to metrized graphs, see the expository article \cite{BF}. The reader should be
aware that metrized graphs also appear in the literature under the names \textit{quantum
graphs}, \textit{metric graphs}, and $c^2$\textit{-networks}.

\section{Definitions, notation, and statement of the main theorem} \label{Notation Section}

A \textit{metrized graph} $\Gamma$ is a compact connected metric space such that for each
point $p \in \Gamma$ there exists a radius $r_p>0$ and a valence $n_p \in \NN$ such that
the open ball in $\Gamma$ of radius $r_p$ about $p$ is isometric to the star-shaped set
$\{re^{2\pi i m/n_p}: 0<r<r_p, 1\leq m \leq n_p\} \subset \CC$ endowed with the path
metric. A \textit{vertex set} $V$ for $\Gamma$ is any finite nonempty subset satisfying
the following properties: (i) $V$ contains all points $p \in \Gamma$ with $n_p \not= 2$;
(ii) for each connected component $U_i \subset \Gamma \setminus V$, the closure $e_i =
\overline{U}_i$ is isometric to a closed interval (not a circle); and (iii) the
intersection $e_i \cap e_j$ consists of at most one point when $i \not= j$.
The set $e_i$ will be called a \textit{segment of $\Gamma$ with respect to the vertex set
$V$}. Given a vertex set $V$ for $\Gamma$, one can associate a combinatorial weighted
graph $G=G(V)$ called a \textit{model} for $\Gamma$. Indeed, index the vertices of $G$ by
the points in $V$ and connect two vertices $p,q$ in $G$ with an edge of weight $1/L$ if
$\Gamma$ has a segment $e$ of length $L$ with endpoints $p,q$. The hypotheses we have
placed on a vertex set ensure that $G$ is a finite connected weighted graph with no
multiple or loop edges.

Given a vertex set $V$ for $\Gamma$, a segment of length $L$ can be isometrically
parametrized by a closed interval $[0,L]$. This parametrization is unique up to a choice
of orientation, and so there is a well-defined notion of Lebesgue measure on the segment
with total mass $L$. Defining the Lebesgue measure for each of the finite number of
segments gives the Lebesgue measure on $\Gamma$, which we denote by $dx$. Evidently it is
independent of the choice of vertex set. For the scope of this paper we will assume that
$\Gamma$ is a fixed metrized graph on which the metric has been scaled so that
$\int_\Gamma dx = 1$; i.e., $\Gamma$ has total length~1.

Choose a signed measure of total mass~1 on $\Gamma$ of the form
\begin{equation} \label{Measure mu} \mu
=\omega(x)dx + \sum_{j=1}^n c_j \delta_{p_j}(x), \end{equation} where $\omega$ is a
real-valued piecewise continuous function in $L^1(\Gamma)$, $c_1, \ldots, c_n$ are real
numbers, and $\delta_{p}(x)$ denotes the unit point mass at $p \in \Gamma$. We may assume
that $X = \{p_1, \ldots, p_n\}$ is a vertex set for $\Gamma$ containing all points where
$\omega$ is discontinuous. This particular vertex set $X$ will be fixed for the rest of
the article. The choice of measure $\mu$ allows some flexibility in applications of the
theory (cf. \S14 of \cite{BR}).

For each $p \in \Gamma$, we define the set $\vect(p)$ of \textit{formal unit vectors
emanating from $p$}. This set has $n_p$ elements in it, where $n_p$ is the valence of
$\Gamma$ at $p$. For $\vec{v} \in \vect(p)$, we write $p + \varepsilon \vec{v}$ for the
point at distance $\varepsilon$ from $p$ in the direction $\vec{v}$, a notion which makes
sense for $\varepsilon$ sufficiently small. Given a function $f:\Gamma \to \CC$, a point
$p \in \Gamma$, and a direction $\vec{v} \in \vect(p)$, the \textit{derivative of $f$ at
$p$ in the direction $\vec{v}$} (or simply \textit{directional derivative}), written
$D_{\vec{v}}f(p)$, is given by
\[
D_{\vec{v}}f(p) = \lim_{\varepsilon\to 0^+} \frac{f(p+\varepsilon\vec{v})
-f(p)}{\varepsilon},
\]
provided this limit exists.

The \textit{Zhang space}, denoted $\Zh(\Gamma)$, is the set of all continuous functions
$f: \Gamma \to \CC$ for which there exists a vertex set $X_f$ for $\Gamma$ such that the
restriction of $f$ to any connected component of $\Gamma \setminus X_f$ is
$\mathcal{C}^2$, and for which $f''(x) \in L^1(\Gamma)$. Directional derivatives of
functions in the Zhang space exist for all points of $\Gamma$ and all directions. The
importance of this space of functions will be made clear in just a moment.

The Laplacian $\Delta(f)$ of a function $f \in \Zh(\Gamma)$ is the \textit{measure}
defined by
\[\Delta(f) = -f''(x)dx - \sum_{p \in \Gamma}\left\{ \sum_{\vec{v} \in \vect(p)}
D_{\vec{v}}f(p)\right\} \delta_p(x).\] One sees that $\sum_{\vec{v} \in \vect(p)}
D_{\vec{v}}f(p) = 0$ for any point $p \in \Gamma \setminus X_f$ so that the outer sum
above is actually finite. For any $f, g \in \Zh(\Gamma)$, the Laplacian operator
satisfies the identities
\begin{equation} \label{Zhang Identity}
\int_\Gamma \overline{g} \ d\Delta(f) = \int_\Gamma f d\overline{\Delta(g)} = \int_\Gamma
f'(x) \overline{g'(x)} dx.\end{equation} Letting $g$ be the constant function with value
$1$ in this last expression shows that the measure $\Delta(f)$ has total mass zero.

Define $\Zh_{\mu}(\Gamma)$ to be the subspace of the Zhang space consisting of all
functions which are orthogonal to the measure $\mu$; i.e., all $f \in \Zh(\Gamma)$ such
that $\int_{\Gamma} f d\mu = 0$. A nonzero function $f \in \Zh_{\mu}(\Gamma)$ is called
an \textit{eigenfunction of the Laplacian (with respect to the measure $\mu$)} if there
exists an \textit{eigenvalue} $\lambda \in \CC$ so that
\[\int_\Gamma  \overline{g} \ d\Delta(f) = \lambda \int_{\Gamma} f(x) \overline{g(x)} dx, \quad
\text{for each $g \in \Zh_{\mu}(\Gamma)$.}\] The integral on the left side of this
equation is called the \textit{Dirichlet inner product} of $f$ and $g$ and is denoted
$\ip{f}{g}_{\Dir}$; the integral on the right side is the usual $L^2$-inner product and
will be denoted $\ip{f}{g}_{L^2}$. Thus we can restate the defining relation for an
eigenfunction as $\ip{f}{g}_{\Dir} = \lambda\ip{f}{g}_{L^2}$. Setting $g=f$ and applying
the relations in \eqref{Zhang Identity} shows that each eigenvalue of the Laplacian
$\Delta$ must be real and positive. The eigenvalues constitute a sequence tending to
infinity; in particular, the dimension of the eigenspace associated to a given eigenvalue
is finite. Let $0<\lambda_1(\Gamma) < \lambda_2(\Gamma) < \lambda_3(\Gamma) < \cdots$
denote the eigenvalues of $\Delta$ with respect to the measure $\mu$. The dimension of
the eigenspace corresponding to $\lambda_i(\Gamma)$ will be denoted $d_i(\Gamma)$. Even
though we have suppressed $\mu$ in the notation for the eigenvalues, they do depend
heavily on the choice of measure.

The Laplacian operator can be defined on a much larger class of continuous functions than
$\Zh(\Gamma)$, and the eigenfunctions of $\Delta$ in general need not lie in the Zhang
space. However, due to the ``smoothness'' of our measure $\mu$, every eigenfunction of
$\Delta$ is indeed $\mathcal{C}^2$ on the complement of our fixed vertex set $X$. In
fact, if $\omega(x)$ (the continuous part of $\mu$) is $\mathcal{C}^m$ on $\Gamma
\setminus X$, then each eigenfunction of the Laplacian is $\mathcal{C}^{m+2}$ on the
complement of $X$. All of this is illuminated by Proposition~15.1 of \cite{BR}.

Now we turn out attention to the discrete approximation of $\Gamma$ by its models. For
each positive integer $N \geq \#X$, choose a vertex set $X \subset V_N \subset \Gamma$ in
such a way that $\#V_N = N$. We will add another hypothesis to this vertex set
momentarily. Define the measure $dx_N$ on $\Gamma$ to be the probability measure with a
point mass of weight $1/N$ at each point in $V_N$. Choose the vertex sets $\{V_N\}$ so
that the sequence of measures $\{dx_N\}$ converges weakly to the Lebesgue measure $dx$.
In doing this, we are modelling the metrized graph $\Gamma$ by finite weighted graphs
$G_N = G(V_N)$ whose vertices become equidistributed in $\Gamma$. This sequence of models
will be of great interest to us.

For each $N$, we also choose a discrete signed measure $\mu_N$ supported on $V_N$ with
total mass~$1$ and finite total variation such that the sequence $\{\mu_N\}$ tends weakly
to $\mu$. For example, we could construct such a sequence in the following way. For each
$p \in V_N$, consider the set of points of $\Gamma$ which are closer to $p$ than to any
other vertex: \[A^N_p=\left\{x \in \Gamma: \dist(x,p)< \dist (x,q) \textrm{ for all
vertices } q \in V_N\setminus \{p\}\right\},\] where $\dist(\cdot, \cdot)$ is the metric
on $\Gamma$. Note that this set is Borel measurable (in fact, it is open). We define the
\textit{discretization} of the measure $\mu$ associated to $V_N$ by setting $\mu_N(p) =
\mu(A^N_p)$ for each $p \in V_N$. The sets $A^N_p$ are pairwise disjoint, so this yields
a discrete signed measure on $\Gamma$ supported on the vertex set $V_N$. As $X \subset
V_N$, we can assert that $\mu_N$ has total mass~1 (because then all point masses of $\mu$
are picked up by some $A^N_p$). Recalling the definition of the measure $\mu$ in
\eqref{Measure mu}, we see that the total variation of $\mu_N$ satisfies
\[|\mu_N|(\Gamma) \leq \int_{\Gamma} |\omega(x)|dx + \sum_{i=1}^n |c_i|.\]
One can check that the measures $\{\mu_N\}$ tend weakly to $\mu$ as $N$ tends to
infinity.

We now wish to tie together the two notions of functions on the vertices of a model $G$
and functions on the metrized graph $\Gamma$. The most natural way to do this is via
linear interpolation of the values of a given function on a vertex set. The class of
\textit{continuous piecewise affine functions} on $\Gamma$, denoted $\CPA(\Gamma)$, is
defined to be the set of all continuous functions $f: \Gamma \to \CC$ for which there
exists a vertex set $X_f$ with the property that $f$ is \textit{affine} on any connected
component of the complement of $X_f$. More precisely, if $U$ is a connected component of
$\Gamma \setminus X_f$, then its closure $e=\overline{U}$ admits an isometric
parametrization $s_e: [0,L] \to e$. To say that $f$ is affine on $e$ is to say that there
are complex constants $A,B$ depending on the segment $e$ so that $f \circ s_e(t) = At+B$.
For each vertex set $V$ of $\Gamma$, we define $\Funct(V)$ to be the subclass of
$\CPA(\Gamma)$ whose values are determined by their values on $V$. That is, $\Funct(V)$
is any continuous function on $\Gamma$ which is a linear interpolation of some set of
values on the vertex set $V$. There is a natural identification of $\Funct(V_N)$ with the
complex vector space $\CC^N$. We also define the $\ell^2$-inner product of two elements
in $\Funct(V_N)$ to be $\ip{f}{g}_{\ell^2} = \int_\Gamma f(x)\overline{g(x)} dx_N$. It is
important to note which value of $N$ we are considering when computing this inner
product.

Each of our preferred models $G_N$ is equipped with a combinatorial weighted Laplacian
matrix (or Kirchhoff matrix), denoted by $Q_N$, which we can view as an abstract linear
operator on the space of functions $\Funct(V_N)$. Label the vertices of $G_N$ by $q_1,
\ldots, q_N$, and assume that edge $q_iq_j=q_jq_i$ has weight $w_{ij}$ (recall that
$w_{ij}$ is the reciprocal of the length of the segment $q_iq_j$ in $\Gamma$ and
$w_{ii}=0$ for all $i$). By definition (cf. \cite{Mo} or \cite{Bo}), the entries of $Q_N$
are given by
\begin{equation} \label{LaplacianMatrixDefinition}
[Q_N]_{ij} = \begin{cases} \sum_k w_{ik}, & \textrm{if $i = j$} \\
-w_{ij}, & \textrm{if $i \not= j$ and $q_i$ is adjacent to $q_j$} \\
0, & \textrm{if $q_i$ is not adjacent to $q_j$.} \end{cases}
\end{equation}
For $f \in \Funct(V_N)$, we write $Q_N f$ or $\{Q_N f\}(x)$ for the unique function in
$\Funct(V_N)$ which has value $\sum_j [Q_N]_{ij} f(q_j)$ at the vertex $q_i$. Evidently
the function $Q_N f$ exhibits the same information one gets from multiplying the matrix
$Q_N$ on the right by the column vector with $j$th entry $f(q_j)$. For functions in
$\Funct(V_N)$, the Laplacian matrix is closely related to the Laplacian operator $\Delta$
via the formula \begin{equation} \label{Laplacian On CPA} \Delta(f) = \sum_{j=1}^n
\{Q_Nf\}(q_j) \ \delta_{q_j}.
\end{equation}
 This is an easy consequence of the definitions of
all of the objects involved (cf. \cite{BF}, Theorem~4).

In analogy with the setup for the continuous Laplacian, we define the class of functions
$\Funct_{\mu_N}(V_N)$ to be the subclass of $\Funct(V_N)$ orthogonal to the measure
$\mu_N$ (a subspace of complex dimension $N-1$). That is, $f \in \Funct_{\mu_N}(V_N)$ if
$\int_{\Gamma} f d\mu_N = 0$. A nonzero function $f \in \Funct_{\mu_N}(V_N)$ will be
called an \textit{eigenfunction for the discrete Laplacian $Q_N$ (with respect to the
measure $\mu_N$}) if there exists an \textit{eigenvalue} $\lambda \in \CC$ so that for
all $g \in \Funct_{\mu_N}(V_N)$, \begin{equation} \label{Discrete Eigendefinition}
\sum_{q \in V_N} \{Q_N f\}(q) \overline{g(q)} = \lambda \sum_{q \in V_N} f(q)
\overline{g(q)}.\end{equation} Using \eqref{Laplacian On CPA}, we can rewrite this
condition as $\ip{f}{g}_{\Dir}= N\lambda \ip{f}{g}_{\ell^2}$, which looks a good deal
like the defining relation for eigenfunctions of $\Delta$. In \S\ref{Integral Operator
section} we show how the eigenfunctions and eigenvalues of $Q_N$ with respect to the
measure $\mu_N$ relate to the usual notions of eigenvalues and eigenvectors. We show that
all of the eigenvalues of $Q_N$ with respect to the measure $\mu_N$ are positive in
Proposition~\ref{positivity}, and in Corollary~\ref{DiscreteBasis} we prove that the
eigenfunctions of $Q_N$ form a basis for $\Funct_{\mu_N}(V_N)$. Let $0< \lambda_{1,N} <
\lambda_{2,N} < \lambda_{3,N} < \cdots$ denote the eigenvalues of $Q_N$, and write
$d_{i,N}$ for the dimension of the eigenspace corresponding to the eigenvalue
$\lambda_{i,N}$.

For each fixed $i \geq 1$ and $N \geq \#X$, let $\HH_N(i) \subset \Funct_{\mu_N}(V_N)$ be
an $\ell^2$-orthonormal basis of the eigenspace of $Q_N$ corresponding to the eigenvalue
$\lambda_{i,N}$. If the eigenspace is empty, set $\HH_N(i) = \emptyset$ (e.g., if $N <
i$). Define $\HH(i) = \bigcup_N \HH_N(i)$. We think of the family $\HH(i)$ as the set of
all eigenfunctions of $Q_N$ corresponding to an $i$th eigenvalue. The family $\HH(i)$ is
obviously not unique, but we fix a choice of $\HH(i)$ for the remainder of the paper.

\begin{thm}[Main Theorem] \label{Main Theorem}
Fix $i \geq 1$. With the hypotheses and conventions as above, we have the following
conclusions:
\begin{enumerate}
\item[(A)]$\lim_{N \to \infty} N \lambda_{i,N} = \lambda_i(\Gamma)$.
\item[(B)] There exists $N_0=N_0\left(i\right)$ so that for all $N>N_0$, the dimension
$d_{i,N}$ of the eigenspace for $Q_N$ with respect to the measure $\mu_N$ corresponding
to the eigenvalue $\lambda_{i,N}$ satisfies $d_{i,N} = d_i(\Gamma)$.
\item[(C)] The family $\HH(i)$ is normal. The subsequential limits of $\HH(i)$ lie in
$\Zh_{\mu}(\Gamma)$ and contain an $L^2$-orthonormal basis for the eigenspace of $\Delta$
corresponding to the eigenvalue $\lambda_i(\Gamma)$.
\end{enumerate}
\end{thm}

As remarked above, we can even sharpen the statement of assertion (C) if the continuous
part of the measure $\mu$ satisfies more smoothness properties. To reiterate, if $\mu =
\omega(x)dx + \sum c_j \delta_{p_j}$ and $\omega$ is $\mathcal{C}^m$ away from the vertex
set $X$, then the subsequential limits of the family $\HH(i)$ are $\mathcal{C}^{m+2}$ on
$\Gamma \setminus X$.

The rate of convergence of eigenvalues and eigenfunctions is not studied in this paper.
Our methods are purely existential, which is unfortunate due to certain interesting
empirical data obtained by the students in the 2003 University of Georgia REU entitled
``Analysis on Metrized Graphs''. The following two conjectures were made under the
hypotheses that the models for $\Gamma$ can be chosen with all edges having equal weights
and that $\mu=dx$:
\begin{itemize}

\item There is a positive constant $M$ such that $ |\lambda_1(\Gamma) -
N\lambda_{1,N} | < M\cdot N^{-3/2}$ for all $N$. (The example in \cite{KoFu1} might lead
one to believe that the error is more like $O(N^{-2})$.)

\item The scaled discrete eigenvalues $N\lambda_{i,N}$ increase \textit{monotonically} to
the limit. This does not appear to be true for more general measures $\mu$.

\end{itemize}
For more information about the REU and other data and conjectures resulting from it, see
the official site: \verb+http://www.math.uga.edu/~mbaker/REU/REU.html+

In the next section we introduce an integral operator whose spectrum is intimately
related to the spectrum of $Q_N$. We use it to show that $\Funct_{\mu_N}(V_N)$ admits a
basis of eigenfunctions of the discrete Laplacian. In \S\ref{Reduction Lemma Section} we
exhibit a reduction of the Main Theorem which is technically simpler to prove. The proof
of the Main Theorem will then proceed by induction on the eigenvalue index $i$ (the case
$i=1$ will be identical to all others). We carry out the induction step in sections
\ref{Normality Section}-\ref{Completion of Proof Section}.

\section{Integral operators and the spectral theory of Laplacians} \label{Integral Operator section}

To begin, we recall the definition of the function $j_{z}(x,y)$. It is given by the
unique continuous solution of $\Delta_x j_{z}(x,y) = \delta_y(x) - \delta_z(x)$ subject
to the initial condition $j_{\zeta}(\zeta,y)=0$ for all $y \in \Gamma$. Here $\Delta_x$
denotes the action of the Laplacian with respect to the variable $x$. As the Laplacian of
the $j$-function has no continuous part, one can show that $j_z(x,y)$ must be piecewise
affine in $x$ for fixed values of $y,z$. It is also true that $j_{z}(x,y)$ is symmetric
in $x$ and $y$, non-negative, jointly continuous in all three variable simultaneously,
and uniformly bounded by~$1$. For elementary proofs of these facts see section~6 of
\cite{BF}.

For the rest of this section, assume $N \geq \#X$ is a fixed integer. We define
$j_{\mu_N}: \Gamma^2 \to \CC$ to be
\begin{equation*}
j_{\mu_N}(x,y) = \int_{\Gamma} j_{z}(x,y) d\mu_N(z) = \sum_{q \in V_N} j_{q}(x,y)
\mu_N(q).
\end{equation*}
Observe that $j_{\mu_N}(x,y)$ is also a piecewise affine function in $x$ for fixed $y$.
We require the following proposition, whose proof is given in \cite{CR} as Lemma~2.16:

\begin{prop} \label{Doesn't Depend on y} Let $\nu$ be a signed Borel measure on $\Gamma$ of finite
total variation. There is a constant $C_{\nu}$ such that for each $y \in \Gamma$,
\[
\int_{\Gamma} j_{\nu}(x,y)d\nu(x) = C_\nu.
\]
\end{prop}

In light of this proposition, we define the kernel function $g_{\nu}:\Gamma^2 \to \CC$ to
be
\begin{equation*}
g_{\nu}(x,y) = j_{\nu}(x,y) - C_\nu.
\end{equation*}
It follows that $\int_{\Gamma} g_{\nu}(x,y) d\nu(y) = 0$. Also, the properties of the
$j$-function mentioned above imply that $g_{\nu}$ is symmetric in its two arguments and
continuous on $\Gamma^2$. Note that since $\Gamma$ is compact, this forces $g_{\nu}$ to
be uniformly continuous on $\Gamma^2$.

Now we recall the fundamental integral transform which effectively inverts the Laplacian
on a metrized graph. For $f \in L^2(\Gamma)$, define
\begin{equation*}
\varphi_{\mu}(f) = \int_\Gamma g_{\mu}(x,y) f(y)dy.
\end{equation*}
It is proved in \cite{BR} that $\varphi_{\mu}$ is a compact Hermitian operator on
$L^2(\Gamma)$ with the property that any element in the image of $\varphi_\mu$ is
orthogonal to the measure $\mu$. In fact, we have the following important equivalence:

\begin{thm}[\cite{BR}, Theorem~12.1] \label{Reciprocal Continuous Eigenvalues}
A nonzero function $f \in \Zh_{\mu}(\Gamma)$ is an eigenfunction of $\varphi_{\mu}$ with
eigenvalue $\alpha > 0$ if and only if $f$ is an eigenfunction of $\Delta$ with respect
to the measure $\mu$ with eigenvalue $\lambda = 1/\alpha > 0$.
\end{thm}

By Proposition~15.1 of \cite{BR} we know that any eigenfunction of $\Delta$ with respect
to the measure $\mu$ lies in $\Zh_{\mu}(\Gamma)$, and the above theorem allows us to
conclude that eigenfunctions of $\varphi_{\mu}$ have the same property.

We now introduce a discrete version of the integral operator whose eigenvalues are
related to the eigenvalues of the Laplacian matrix $Q_N$ in much the same way as in the
above theorem. Define the \textit{discrete integral operator} $\p:\Funct(V_N) \to
\Funct(V_N)$ by
\[\p(h)(x) = \int_\Gamma g_{\mu_N}(x,y) h(y) dy_N, \quad \textrm{where $x \in V_N$.}\]
The defining equation for $\p(h)$ only gives its values at the vertices, but recall that
any function in $\Funct(V_N)$ is determined by its values on $V_N$ by linear
interpolation.

The next lemma shows that the Laplacian matrix is essentially a left inverse for $\p$, up
to a scaling factor and a correction term.

\begin{prop} \label{Laplacian Inverse}
If $f \in \Funct(V_N)$, then for any $q \in V_N$, we have
\begin{equation*}
\left\{Q_N \p\left(f\right)\right\}(q) = \frac{1}{N}f(q) - \left(\int_{\Gamma}f(x) \
dx_N\right)\mu_N(q).
\end{equation*}
\end{prop}

\begin{proof} This is a restatement of Proposition~7.1 of \cite{BR} (setting $\nu = f(x)dx_N$,
$\mu=\mu_N$, and using equation~\eqref{Laplacian On CPA} to relate the Laplacian $\Delta$
to the discrete Laplacian matrix $Q_N$).
\end{proof}

Now we exhibit two useful facts about the kernel and image of the operator $\p$.

\begin{lem} \label{Trivial Kernel}
$\Ker(\p) \bigcap \Funct_{\mu_N}(V_N) = 0$
\end{lem}

\begin{proof}
Suppose $f \in \Ker(\p) \bigcap \Funct_{\mu_N}(V_N)$.  By Proposition~\ref{Laplacian
Inverse}, we find for each $q \in V_N$ that
\begin{equation*} \label{Kernel Relation}
0 = \left\{Q_N\p\left(f\right)\right\}(q) = \frac{1}{N}f(q) - \left(\int_{\Gamma}f(x) \
dx_N\right)\mu_N(q).
\end{equation*}
Put $C_f=\int_{\Gamma}f(x) dx_N$. If $C_f=0$, then $f \equiv 0$. If $C_f \not= 0$, then
the above equality implies $f(q) = NC_{f}\mu_N(q)$. But we then obtain the following
contradiction to the fact that $f \in \Funct_{\mu_N}(V_N)$:
\begin{equation*}
\int_{\Gamma} f(x) d\mu_N(x) = NC_{f} \sum_{q_i \in V_N} \mu_N(q_i)^2 \not = 0.
\end{equation*}
\end{proof}

We remark that the previous proof actually shows the kernel of $\p$ acting on the space
$\Funct(V_N)$ consists of all scalar multiples of the piecewise affine function defined
by $q \mapsto \mu_N(q), q \in V_N$.

\begin{lem} \label{Right Space} The operator $\p$ is $\ell^2$-Hermitian with image in $\Funct_{\mu_N}(V_N)$.
That is, if $f \in \Funct(V_N)$, then
\begin{equation*}
\int_{\Gamma} \p\left(f\right)(x) d\mu_N(x) = 0.
\end{equation*}
\end{lem}

\begin{proof} Suppose $f,g \in \Funct(V_N)$. As $g_{\mu_N}$ is real and symmetric, we see
\begin{equation*}
\begin{aligned}
\ip{\p(f)}{g}_{\ell^2} &= \int_\Gamma \left( \int_\Gamma g_{\mu_N} (x,y) f(y) dy_N\right)
\overline{g(x)} dx_N \\
&= \int_\Gamma f(y) \overline{\left(\int_{\Gamma} g_{\mu_N}(x,y) g(x) dx_N\right)} dy_N
\\
&= \ip{f}{\p(g)}_{\ell^2}.
\end{aligned}
\end{equation*} Now recall that for any fixed $y \in \Gamma$, the function $g_{\mu_N}(x,y)$ is orthogonal
to the measure $\mu_N$. Thus
\begin{equation*}
\begin{aligned}
\int_\Gamma \p(f)(x) d\mu_N(x) &=\int_\Gamma \left( \int_\Gamma g_{\mu_N} (x,y) f(y)
dy_N\right)
d\mu_N(x) \\
&= \int_\Gamma f(y) \left(\int_{\Gamma} g_{\mu_N}(x,y)  d\mu_N(x)\right) dy_N = 0.
\end{aligned}
\end{equation*}
\end{proof}

It is interesting to pause for a moment to see how the eigenvalues and eigenfunctions of
$Q_N$ with respect to the measure $\mu_N$ relate to the usual notions of eigenvalues and
eigenvectors.

\begin{prop} \label{Equivalent Laplacian}
The function $f \in \Funct_{\mu_N}(V_N)$ is an eigenfunction for $Q_N$ with respect to
the measure $\mu_N$ if and only if there exists a constant $\lambda \in \CC$ such that
for all vertices $q \in V_N$,
\begin{equation} \label{First Claim}
\{Q_N f\}(q) = \lambda \left\{f(q) - N\left(\int_\Gamma f(x) dx_N\right)
\mu_N(q)\right\}.
\end{equation}
In particular, $f$ is an eigenfunction for $Q_N$ with respect to the measure $dx_N$ if
and only if $f$ is a non-constant eigenfunction for $Q_N$ in the usual sense of a linear
operator.
\end{prop}

\begin{proof}
The final claim follows from \eqref{First Claim} because the integral term vanishes from
the result when $f$ is orthogonal to the measure $dx_N$.

If $f \in \Funct_{\mu_N}(V_N)$ is an eigenfunction for $Q_N$ with respect to the measure
$\mu_N$, then there is $\lambda \in \CC$ such that for all $g \in \Funct_{\mu_N}(V_N)$,
\[\sum_{q \in V_N} \left\{Q_N f\right\}(q) \overline{g(q)} = \lambda \sum_{q \in V_N} f(q)
\overline{g(q)}.\] We can rewrite this relation as $N\ip{Q_Nf - \lambda f}{g}_{\ell^2} =
0$. Setting $F=Q_Nf-\lambda f$, we conclude that $F$ is in the $\ell^2$-orthogonal
complement of the space $\Funct_{\mu_N}(V_N)$ (as a subspace of $\Funct(V_N)$). One
easily sees that the function $q \mapsto \mu_N(q)$ for $q \in V_N$ is a basis of the
orthogonal complement. Thus $F(q) = M\mu_N(q)$ for some constant $M$ and all $q \in V_N$.

The value $\{Q_Nf\}(q)$ can be interpreted as the weight of the point mass of $\Delta(f)$
at the point $q$ using equation~\eqref{Laplacian On CPA}. As the Laplacian is always a
measure of total mass zero, and the measure $\mu_N$ has total mass~1, we may sum the
equation $F(q)=M\mu_N(q)$ over all vertices $q$ to see that \[M=-\lambda \sum_{q\in V_N}
f(q) = -N\lambda \int_\Gamma f(x) dx_N.\] This finishes the proof in one direction. The
other direction is an immediate computation.
\end{proof}

\begin{prop} \label{positivity} The eigenvalues of $Q_N$ acting
on $\Funct(V_N)$ are nonnegative. The kernel of $Q_N$ is $1$-dimensional with basis the
constant function with value~$1$. The eigenvalues of $Q_N$ with respect to the measure
$\mu_N$ are all positive.
\end{prop}

\begin{proof}
Using equation~\eqref{Laplacian On CPA} and the notation for the weights on the edges of
the model $G_N$ in \eqref{LaplacianMatrixDefinition}, we can expand the Dirichlet norm as
\begin{equation*}
\begin{aligned}
\|f_N\|^2_{\Dir} &= \int_{\Gamma} \overline{f_N(x)} \  d\Delta(f_N)
= \int_\Gamma \overline{f_N(x)} \left[\sum_{i=1}^N \left\{ Q_Nf\right\}(q_i) \ \delta_{q_i} (x)\right] \\
&= \sum_{q_i \in V_N}  \overline{f_N(q_i)} \left\{Q_Nf_N \right\}(q_i)  \\ &= \sum_{q_i
\in V_N} \left( \left|f_N(q_i)\right|^2 \sum_{q_k \in V_N} w_{ik} - \overline{f_N(q_i)}
\sum_{q_j \in V_N} f_N(q_j)w_{ij} \right).
\end{aligned}
\end{equation*}
Note that $w_{ii}=0$ since $G_N$ has no loop edges. Let us rearrange this last sum to be
over the edges of $G_N$ instead of over its vertices. Observe that summing over all pairs
of vertices as above is equivalent to summing twice over all edges of the graph. We count
$w_{ik}|f_N(q_i)|^2$ for each end of an edge $q_iq_k$. We also count
$-w_{ij}f_N(q_i)\overline{f_N(q_j)}$ and $-w_{ij}\overline{f_N(q_i)}f_N(q_j)$ for each
edge $q_iq_j$. Here we have implicitly taken advantage of the symmetry of $Q_N$. This
yields
\begin{equation}
\begin{aligned}
\label{Laplacian Identity} \|f_N\|^2_{\Dir} &= \sum_{\textrm{edges } q_iq_k } w_{ik}
\left(\left|f_N(q_i)\right|^2  + \left|f_N(q_k)\right|^2\right) \\ & \hspace{0.6in}-
\sum_{\textrm{edges }
q_iq_j } w_{ij} \left(f_N(q_i)\overline{f_N(q_j)} +\overline{f_N(q_i)}f_N(q_j) \right)\\
&= \sum_{\textrm{edges } q_iq_j} w_{ij} \left| f_N(q_i) -f_N(q_j) \right|^2.
\end{aligned}
\end{equation}
Note that in these sums we count edge $q_iq_j=q_jq_i$ only once.

If $f \in \Funct(V_N)$ is an $\ell^2$-normalized eigenfunction for $Q_N$ with respect to
the measure $\mu_N$ with eigenvalue $\lambda$, then equation~\eqref{Laplacian Identity}
and the defining relation for an eigenfunction shows that $0 \leq \|f\|_{\Dir}^2 =
N\lambda$. Thus, $\lambda$ is nonnegative. If $\lambda=0$, then a more careful look at
\eqref{Laplacian Identity} shows that our eigenfunction $f$ satisfies $f(q_i)=f(q_j)$ for
every pair of adjacent vertices $q_i, q_j$. Since $G_N$ is connected, we see that $f$ is
constant on the vertex set $V_N$. This proves the second assertion.

Suppose further that $f$ has constant value $M$. Then $\int_{\Gamma} f d\mu_N = M$. If $f
\in \Funct_{\mu_N}(V_N)$, we are forced to conclude that $M=0$. That is, our function $f$
is identically zero on $V_N$. This is horribly contrary to our definition of an
eigenfunction, and so we conclude that $\lambda > 0$ in this case.
\end{proof}

A function $f \in \Funct(V_N)$ is called an \textit{eigenfunction} for the operator $\p$
if there exists an \textit{eigenvalue} $\alpha \in \CC$ such that $\p(f) = \alpha f$.
Lemma~\ref{Right Space} implies that any eigenfunction with nonzero eigenvalue must lie
in $\Funct_{\mu_N}(V_N)$. Now we prove the theorem which relates the nonzero eigenvalues
and eigenfunctions of the integral operator $\p$ to the discrete Laplacian $Q_N$.

\begin{thm} \label{reciprocate}
A nonzero function $f \in \Funct_{\mu_N}(V_N)$ is an eigenfunction of $Q_N$ (with respect
to the measure $\mu_N$) with eigenvalue $\lambda$ if and only if it is an eigenfunction
of $\p$ with eigenvalue $\frac{1}{N\lambda}$.
\end{thm}

\begin{proof}
Suppose $f$ is an eigenfunction of $\p$ in $\Funct_{\mu_N}(V_N)$ with eigenvalue
$\alpha$.  Proposition~\ref{Laplacian Inverse} implies that for each $q \in V_N$
\begin{equation*}
\alpha \{Q_Nf\}(q)  = \left\{Q_N\p(f)\right\}(q) = \frac{1}{N}f(q) - \left(\int_{\Gamma}
f(x)dx_N\right) \mu_N(q).
\end{equation*}
For $g \in \Funct_{\mu_N}(V_N)$, we may multiply this last equality by $\overline{g(q})$
and sum over all vertices $q \in V_N$ to see that
\begin{equation*}
\alpha\sum_{q \in V_N}\{Q_N f\}(q)\overline{g(q)} = \frac{1}{N}\sum_{q \in V_N}
f(q)\overline{g(q)}.
\end{equation*}
The last sum vanishes because $g$ is orthogonal the the measure $\mu_N$. If $\alpha=0$,
then $f\equiv 0$ by Lemma~\ref{Trivial Kernel}; but eigenfunctions are not identically
zero. Thus $f$ satisfies the defining equation~\eqref{Discrete Eigendefinition} for an
eigenfunction of $Q_N$ with eigenvalue $\frac{1}{N\alpha}$.

Conversely, let $f \in \Funct_{\mu_N}(V_N)$ be an eigenfunction of $Q_N$ with eigenvalue
$\lambda$. Multiplying the result of Proposition~\ref{Laplacian Inverse} by $N\lambda$
and subtracting from the result in Proposition~\ref{Equivalent Laplacian} shows
\[Q_N\left\{f-N\lambda \p(f)\right\}  \equiv 0.\]
By Proposition~\ref{positivity}, the kernel of $Q_N$ is precisely the constant functions
on $\Gamma$. Hence $f-N\lambda \p(f) \equiv M$ for some constant $M$. As $f$ and $\p(f)$
are orthogonal to the measure $\mu_N$, we know that integrating the equation $f-N\lambda
\p(f) = M$ against $\mu_N$ produces $M=0$. Since all eigenvalues of $Q_N$ with respect to
the measure $\mu_N$ are nonzero, we deduce that $\p(f) = \frac{1}{N\lambda}f$.
\end{proof}

\begin{cor}
All of the eigenvalues of $\p$ acting on $\Funct_{\mu_N}(V_N)$ are positive.
\end{cor}

\begin{proof} Apply the previous theorem and Proposition~\ref{positivity}.
\end{proof}

\begin{cor} \label{DiscreteBasis}There exists a basis of $\Funct_{\mu_N}(V_N)$ consisting of eigenfunctions of
$Q_N$ with respect to the measure $\mu_N$.
\end{cor}

\begin{proof} The space $\Funct_{\mu_N}(V_N)$ admits a basis of eigenfunctions of $\p$ by
Lemmas~\ref{Right Space}~and~\ref{Trivial Kernel} and the finite dimensional spectral
theorem. Now apply the previous theorem.
\end{proof}

\section{Reduction to a weaker form of the main theorem} \label{Reduction Lemma Section}

Before embarking on the proof of the Main Theorem, we show the following reduction:

\begin{lem}[Reduction Lemma]  \label{reduction lemma} Suppose that for each fixed $i \geq 1$ the following
assertions are true:
\begin{enumerate}
\item[(A)]$\lim_{N \to \infty} N \lambda_{i,N} = \lambda_i(\Gamma)$.
\item[(B')] There exists $N_0=N_0\left(i\right)$ so that for all $N>N_0$, the dimension
of the eigenspace for $Q_N$ corresponding to the eigenvalue $\lambda_{i,N}$ satisfies
$d_{i,N} \leq d_i(\Gamma)$.
\item[(C')] The family $\HH(i)$ is normal. The subsequential limits of $\HH(i)$ lie in
$\Zh_{\mu}(\Gamma)$, have unit $L^2$-norm, and are eigenfunctions of $\Delta$
corresponding to the eigenvalue $\lambda_i(\Gamma)$.
\end{enumerate}
Then Theorem~\ref{Main Theorem} is true.
\end{lem}

The proof of the Reduction Lemma will require a few preliminary results, which will take
up the majority of this section.

\begin{lem} \label{Dirichlet Convergence Lemma} Suppose $f \in \Zh(\Gamma)$ with the property
that $f$ is $\mathcal{C}^2$ on $\Gamma \setminus X$. Define $f_N$ to be the unique
function in $\Funct(V_N)$ which agrees with $f$ on the vertex set $V_N$. Then
\begin{equation*}
 \lim_{N \to \infty} \|f_N\|_{\Dir} = \|f\|_{\Dir}.
 \end{equation*}

\end{lem}

\begin{proof} We begin by recalling equation~\eqref{Laplacian Identity}:
\begin{equation} \label{Laplacian Identity2}
\|f_N\|^2_{\Dir} = \sum_{\textrm{edges } q_iq_j} w_{ij} \left| f_N(q_i) -f_N(q_j)
\right|^2.
\end{equation}
As $\Gamma$ can be decomposed into a finite collection of segments $\{e\}$using the
preferred vertex set $X$, we can sum over the segments of $\Gamma$ to get
\begin{equation}
 \label{slick sum}
\|f_N\|^2_{\Dir} = \sum_{\substack{\textrm{segments $e$} \\ \textrm{of $\Gamma$}}}
\hspace{.1 in} \sum_{\substack{\textrm{edges } q_iq_j \\ \textrm{of $G_N$ with $q_i,q_j
\in e$}}} w_{ij} \left| f_N(q_i) -f_N(q_j) \right|^2.
\end{equation}

The integral is linear, and there are only finitely many segments of $\Gamma$ over which
we wish to integrate, so it suffices to show that the inner sum in \eqref{slick sum}
converges to $\int_e |f'(x)|^2 dx$ for any segment $e$ of $\Gamma$ (by
equation~\eqref{Zhang Identity}). Hence, we may assume for the remainder of this proof
that $\Gamma$ consists of exactly one segment $e$; that is, $\Gamma$ is isometric to a
closed interval of length~1.

The single segment $e$ of $\Gamma$ admits an isometric parametrization $s_e:[0,1] \to e$.
The vertex set $V_N$ corresponds to a partition $0=t_1 < t_2 < \cdots < t_N = 1$. We
write $f_e = f \circ s_e$ for ease of  notation. Now \eqref{Laplacian Identity2} takes
the form
\begin{align*} \label{monotonic equality}
\|f_N\|^2_{\Dir} &= \sum_{i=1}^{N-1} w_{i(i+1)}
\left|f_e(t_{i+1}) - f_e(t_i)\right|^2 \\
&= \sum_{i=1}^{N-1} \left(t_{i+1} -
t_i\right)\left|\frac{f_e(t_{i+1})-f_e(t_i)}{t_{i+1}-t_i}\right|^2.
\end{align*}

For each $i$ there is $t_i^* \in (t_i,t_{i+1})$ so that $f_e(t_{i+1})-f_e(t_i) =
f_e'(t_i^*) (t_{i+1}-t_i)$ by the mean value theorem. Hence
\begin{align*}
\|f_N\|^2_{\Dir} &= \sum_{i=1}^{N-1} \left(t_{i+1} - t_i\right)
\left|f_e'(t_i^*)\right|^2.
\end{align*}
 But this last expression is just
a Riemann approximation to $\int_0^1 |f_e'(x)|^2 dx$. By Proposition~5.2 in \cite{BR}, we
find that $f_e'$ is a continuous function on $[0,1]$ so that as $N$ tends to infinity
these approximations actually do limit to the desired integral.

\end{proof}

The following immediate corollary won't be of any use to us in our present task, but it
is nice to include for completeness.

\begin{cor} \label{Dirichlet Inner Product Convergence}If $f,g \in \Zh(\Gamma)$, then we have the limit
\begin{equation*} \lim_{N \to \infty} \ip{f_N}{g_N}_{\Dir} \longrightarrow \ip{ f}{g}_{\Dir},
\end{equation*}
where $f_N$ and $g_N$ are affine approximations of $f$ and $g$ as in the statement of
Lemma~\ref{Dirichlet Convergence Lemma}.
\end{cor}

\begin{proof}
This is an easy consequence of Lemma~\ref{Dirichlet Convergence Lemma} and the
polarization identity
\[\ip{f_N}{g_N}_{\Dir} = \frac{1}{4} \sum_{n=0}^3 i^n \|f_N + i^n g_N\|^2_{\Dir},\]
which may be checked easily by expanding the norms on the right-hand side. Here $i$
denotes a fixed complex root of $-1$.
\end{proof}

Here is a technical lemma that will be needed in the Approximation Lemma
(Lemma~\ref{Precision Lemma}) and again in later sections.

\begin{lem} \label{Technical Crap 5}
Suppose for each $N$ that $k_N$ is a function in $\Funct(V_N)$ such that the sequence
$\{k_N\}$ converges uniformly to a (continuous) function $k: \Gamma \to \CC$. Then
$\|k_N\|_{\ell^2} \to \|k\|_{L^2}$ as $N \to \infty$. If $\{k'_N\}$ is another such
sequence of functions converging uniformly to some $k': \Gamma \to \CC$, then
$\ip{k_N}{k'_N}_{\ell^2} \to \ip{k}{k'}_{L^2}$ as $N \to \infty$.
\end{lem}

\begin{proof}
Suppose $N$ is so large that $k_N$ is uniformly within $\varepsilon$ of $k$. On one hand,
we have
\begin{equation*}
\begin{aligned}
\|k_N\|_{\ell^2}^2 &= \int_{\Gamma} |k_N(x)|^2 dx_N \\
&\leq \int_{\Gamma} \left(|k(x)| + \varepsilon \right)^2 dx_N
\\
&= \|k\|_{\ell^2}^2 + 2\varepsilon\|k\|_{\ell^1} + \varepsilon^2 \\
&\longrightarrow \|k\|_{L^2}^2 + 2\varepsilon\|k\|_{L^1} + \varepsilon^2.
\end{aligned}
\end{equation*}
The convergence in the final step follows from weak convergence of $dx_N$ to $dx$.
Letting $\varepsilon \to 0$ shows that $\limsup_{N \to \infty} \|k_N\|_{\ell^2} \leq
\|k\|_{L^2}$.

Similarly,
\begin{equation*}
\begin{aligned}
\int_{\Gamma} |k_N(x)|^2 dx_N &\geq \int_{\Gamma} \left(|k(x)| - \varepsilon
\right)^2 dx_N \\
&= \|k\|_{\ell^2}^2 - 2\varepsilon\|k\|_{\ell^1} + \varepsilon^2 \\
&\longrightarrow \|k\|_{L^2}^2 - 2\varepsilon\|k\|_{L^1} + \varepsilon^2.
\end{aligned}
\end{equation*}
Now we see that $\|k\|_{L^2} \leq \liminf_{N \to \infty} \|k_N\|_{\ell^2}$.

The final statement follows by the polarization identity as in the proof of
Corollary~\ref{Dirichlet Inner Product Convergence}.
\end{proof}

\begin{lem}[Approximation Lemma] \label{Precision Lemma} Fix $m \geq 1$ and
suppose $f \in \Zh_{\mu}(\Gamma)$ is an $L^2$-normalized eigenfunction of $\Delta$ with
corresponding eigenvalue $\lambda_m(\Gamma)$. For each $N$, let $f_N$ be the unique
function in $\Funct(V_N)$ which agrees with $f$ at the vertices in $V_N$. Assume also
that assertions (A), (B'), and (C') of the Reduction Lemma (Lemma~\ref{reduction lemma})
are satisfied for $j \leq m-1$. Given any subsequence of $\{G_N\}$, there exists a
further subsequence such that for each $\varepsilon> 0$ there is a positive integer $N_1$
with the property that for every $N>N_1$ with $G_N$ in our sub-subsequence, we can find a
function $\widetilde{f}_N \in \Funct_{\mu_N}(V_N)$ for which the following conditions
hold simultaneously:
\begin{enumerate}
\item[(i)] $\|\widetilde{f}_N\|_{\ell^2} = 1$;
\item[(ii)] $\widetilde{f}_N$ is $\ell^2$-orthogonal to every $h \in \HH_N(j)$, $j=1, \ldots,
m-1$; and
\item[(iii)] $\left| \|\widetilde{f}_N\|_{\Dir}^2 - \|f_N\|_{\Dir}^2
\right| < \varepsilon$.
\end{enumerate}
\end{lem}

\begin{proof}
As condition (B') of the Reduction Lemma holds, and each $d_j(\Gamma)$ is finite, we may
pass to a subsequence of models $\{G_N\}$ such that $d_{j,N}$ is independent of $N$ for
every $j \leq m-1$ and all $N$ sufficiently large in the subsequence. This implies that
for large $N$, the set $\TT_N=\bigcup_{j=1}^{m-1} \HH_N(j)$ is finite with cardinality
independent of $N$. Denote this cardinality by $T$, and let us label the elements of
$\TT_N$ as $h_N^1, \ldots, h_N^T$. The set $\TT_N$ is an $\ell^2$-orthonormal system by
definition of the sets $\HH_N(j)$ and the fact that eigenfunctions corresponding to
distinct eigenvalues are $\ell^2$-orthogonal (a consequence of the definition of
eigenfunction).

Write $B_N = \int_{\Gamma} f d\mu_N$; we will abuse notation in what follows and let
$B_N$ also denote the constant function on $\Gamma$ with value $B_N$. Define
\begin{equation*} \label{Tilde Guy}
\widetilde{f}_N = \frac{f_N - \sum_{j=1}^T \ip{f_N}{h_N^j}_{\ell^2}h_N^j -
B_N}{\left\|f_N - \sum_{j=1}^T \ip{f_N}{h_N^j}_{\ell^2}h_N^j - B_N\right\|_{\ell^2}}.
\end{equation*}
It will be shown in a moment that for $N$ sufficiently large lying in a certain
subsequence, the denominator of this expression is nonzero.

The definition of $\HH(j)$ implies that each $h_N^j$ is orthogonal to the measure
$\mu_N$, so a small calculation shows $\widetilde{f}_N$ lies in $\Funct_{\mu_N}(V_N)$.
Evidently condition (i) holds. Also, we note that $B_N$ is an eigenfunction of $\p$
acting on $\Funct(V_N)$ with corresponding eigenvalue zero. Any two eigenfunctions of a
self-adjoint operator that have distinct associated eigenvalues are orthogonal, and hence
$B_N$ must be $\ell^2$-orthogonal to each $h_N^j$. Another quick calculation shows that
$h_N^j$ is orthogonal to $\widetilde{f}_N$ for all $j \leq m-1$, which is condition (ii).

It remains to check that condition (iii) holds for our choice of $\widetilde{f}_N$ upon
passage to a further subsequence. First observe that the definition of an eigenfunction
for $Q_N$ implies that distinct elements of $\TT_N$ are orthogonal with respect to the
Dirichlet inner product. Also, $\ip{B_N}{g}_{\Dir} = 0$ for all $g \in \Funct(V_N)$ since
constant functions are in the kernel of $Q_N$. Expanding the Dirichlet norm and
simplifying using these two observations gives
\begin{equation*}
\begin{aligned}
\left\|f_N - \sum_{j=1}^T \ip{f_N}{h_N^j}_{\ell^2}h_N^j - B_N\right\|_{\Dir}^2 =
\left\|f_N\right\|^2_{\Dir} - \sum_{j=1}^T
N\zeta_{j,N}\left|\ip{f_N}{h_N^j}_{\ell^2}\right|^2,
\end{aligned}
\end{equation*}
where $\zeta_{j,N}$ is the eigenvalue associated to the eigenfunction $h_N^j$. Of course,
$\zeta_{j,N} = \lambda_{k,N}$ for some $k \leq m-1$.

We have $T$ sequences of functions $\{h_N^1\}, \ldots, \{h_N^T\}$, each of which lies in
a particular $\HH(j)$. The normality of each $\HH(j)$ allows us to find a further
subsequence of models and a set of functions $\TT=\left\{h^1, \ldots, h^T\right\}$ so
that $h_N^j \to h^j$ uniformly on $\Gamma$ along this subsequence. For the rest of the
proof we will assume that $N$ lies in the subsequence we have just specified.

As $\TT_N$ forms an $\ell^2$-orthonormal set for each $N$, we see that $\TT$ forms an
$L^2$-orthonormal set of eigenfunctions for $\Delta$ (assertion (C')). The associated
eigenvalues of these functions are all strictly smaller than $\lambda_m(\Gamma)$, and
hence each $h^j$ is $L^2$-orthogonal to our initial function $f$. Lemma~\ref{Technical
Crap 5} implies that the inner product $\ip{f_N}{h_N^j}_{\ell^2}$ tends to
$\ip{f}{h^j}_{L^2}=0$ as $N$ tends to infinity.

By Lemma~\ref{Dirichlet Convergence Lemma}, weak convergence of $\{dx_N\}$ to $dx$, the
previous paragraph, and the hypothesis that the scaled eigenvalues $N\lambda_{j,N}$
converge to $\lambda_j(\Gamma)$ for each $j \leq m-1$, we conclude that
\begin{equation} \label{Numerator estimate}
\left\|f_N\right\|^2_{\Dir} - \sum_{j=1}^T
N\zeta_{j,N}\left|\ip{f_N}{h_N^j}_{\ell^2}\right|^2 \longrightarrow \|f\|^2_{\Dir} -
\sum_{j=1}^T \zeta_j\left|\ip{f}{h^j}_{L^2}\right|^2 =\|f\|^2_{\Dir}
\end{equation}
as $N$ tends to infinity.

We also note that the function $\sum_{j=1}^T\ip{f_N}{h_N^j}_{\ell^2}h_N^j + B_N$ tends to
zero uniformly. Here $B_N \to 0$ by weak convergence of $\{\mu_N\}$ to $\mu$. The
Minkowski inequality and Lemma~\ref{Technical Crap 5} imply that
\begin{equation*}
\begin{aligned}
\left\|f_N - \sum_{j=1}^T\ip{f_N}{h_N^j}_{\ell^2}h_N^j - B_N\right\|_{\ell^2} &\leq
\left\|f_N\right\|_{\ell^2} + \left\|\sum_{j=1}^T\ip{f_N}{h_N^j}_{\ell^2}h_N^j +
B_N\right\|_{\ell^2} \\
&= \|f\|_{L^2} + o(1),
\end{aligned}
\end{equation*}
where the error term depends only on $N$. The Minkowski inequality also gives an
identical lower bound so that
\begin{equation} \label{Denominator estimate}
\left\|f_N - \sum_{j=1}^T\ip{f_N}{h_N^j}_{\ell^2}h_N^j - B_N\right\|_{\ell^2} \to
\|f\|_{L^2} = 1
\end{equation} as $N$ goes to infinity. This shows that the denominator of
$\widetilde{f}_N$ is in fact nonzero for $N$ sufficiently large.

Finally, we use the convergence relations \eqref{Numerator estimate} and
\eqref{Denominator estimate} along with Lemma~\ref{Dirichlet Convergence Lemma} to
conclude that as $N$ tends to infinity
\begin{equation*}
 \left\|\widetilde{f}_N\right\|_{\Dir}^2 -
\left\|f_N\right\|_{\Dir}^2  = \frac{\left\|f_N\right\|^2_{\Dir} - \sum_{j=1}^T
N\zeta_{j,N}\left|\ip{f_N}{h_N^j}_{\ell^2}\right|^2}{\left\|f_N -
\sum_{j=1}^T\ip{f_N}{h_N^j}_{\ell^2}h_N^j - B_N\right\|_{\ell^2}^2} - \|f_N\|^2_{\Dir}
\longrightarrow 0.
\end{equation*}
This completes the proof of condition (iii).
\end{proof}

And now we require one final proposition which gives us a characterization of the $m$th
largest eigenvalue of $Q_N$ as a minimum of the Dirichlet norm on the $\ell^2$ unit
circle.

\begin{lem} \label{Min at Eigenvector} Fix $N \geq \#X$.
For a given $m \geq 1$, let $\TT_N = \bigcup_{j=1}^{m-1} \HH_N(j)$. Denote by
$\TT_N^{\perp}$ the $\ell^2$-orthogonal complement of the span of $\TT_N$ inside the
space $\Funct_{\mu_N}(V_N)$. Provided that $\TT_N^{\perp} \not= \emptyset$, we have
\begin{equation*}
N\lambda_{m,N} = \min_{\substack{g \in \TT_N^{\perp} \\
g \not\equiv 0}} \frac{\|g\|_{\Dir}^2}{\|g\|_{\ell^2}^2}
\end{equation*}
\end{lem}

\begin{proof} As $\Funct_{\mu_N}(V_N)$ admits a basis of
eigenfunctions of $Q_N$ (Corollary~\ref{DiscreteBasis}), there exists an
$\ell^2$-orthonormal basis $\{f_1, \ldots, f_r\}$ for $\TT_N^{\perp}$. Let us assume that
the corresponding eigenvalues are $\gamma_1 \leq \gamma_2 \leq \cdots \leq \gamma_r$.
Note that $\gamma_1=\lambda_{m,N}$ since $\TT_N^{\perp}$ consists of all of the
eigenfunctions of $Q_N$ with eigenvalue strictly larger than $\lambda_{m-1,N}$.

To prove the proposition it suffices to consider the minimum over all $g \in
\TT_N^{\perp}$ with unit $\ell^2$-norm. Using our basis, write $g = \sum_{j=1}^r a_j
f_j$, where the complex numbers $a_j$ satisfy $\sum_j |a_j|^2 = 1$. Since
$\ip{f_j}{f_j}_{\Dir} = N \gamma_j$, we see that
\begin{equation*}
\ip{g}{g}_{\Dir} = \sum_{j=1}^r |a_j|^2 \ip{f_j}{f_j}_{\Dir} = N\sum_{j=1}^r |a_j|^2
\gamma_j \geq N \gamma_1\sum_{j=1}^r |a_j|^2 = N \lambda_{m,N}.
\end{equation*}
Note that if we take $g \in \Funct_{\mu_N}(V_N)$ to be an $\ell^2$-normalized
eigenfunction of $Q_N$ with eigenvalue $\lambda_{m,N}$, then equality is actually
achieved in this last computation.
\end{proof}

Finally, we may now return to the
\begin{proof}[Proof of Reduction Lemma]
We assume that conditions (A), (B'), and (C') hold for every $i \geq 1$. Clearly we only
need to check that assertions (B) and (C) in the Main Theorem are true. If assertion (B)
fails for some value of $i$, then we may choose a minimum such $i$. There exists a
positive integer $D$ strictly smaller than $d_i(\Gamma)$ and an infinite subsequence of
models $\{G_N\}$ with $d_{i,N}=D$ for all $G_N$ in the subsequence. For each such $N$,
select distinct functions $h_N^1, \ldots, h_N^{D} \in \HH_N(i)$. Using the normality of
$\HH(i)$ as given by hypothesis (C'), we may pass to a further subsequence and assume
that there are $D$ limit functions $h^1, \ldots, h^{D}$, with $h_N^j \to h^j$ uniformly
along this subsequence for each $j$. Again using assertion (C') we find that each $h^j$
is an eigenfunction of $\Delta$ with associated eigenvalue $\lambda_i(\Gamma)$. As $D <
d_i(\Gamma)$, there is some eigenfunction of $\Delta$ with eigenvalue $\lambda_i(\Gamma)$
which is not in the span of $\{h^1, \ldots, h^{D}\}$. Choose such a function with unit
$L^2$-norm and denote it by $f$.

Now choose a sub-subsequence for which the conclusion of the Approximation Lemma is
satisfied for the function $f$ in the case $m=i+1$. Given $\varepsilon > 0$ and $N$
large, we pick $\widetilde{f}_N$ as in the Approximation Lemma. We apply Lemma~\ref{Min
at Eigenvector} in the case $m=i+1$ to see that
\begin{equation*}
N\lambda_{i+1,N} = \min_{\substack{g \in \TT_N^{\perp} \\
g \not\equiv 0}} \frac{\|g\|_{\Dir}^2}{\|g\|_{\ell^2}^2} \leq
\|\widetilde{f}_N\|^2_{\Dir} < \|f_N\|^2_{\Dir}+  \varepsilon.
\end{equation*}
We remark that $\TT_N^{\perp} \not= \emptyset$ for $N$ sufficiently large since its
dimension is at least $N-1 - \sum_{j=1}^{i} d_j(\Gamma)$ by assertion (B'). Letting $N$
go to infinity through our subsequence and using assertion (A) and Lemma~\ref{Dirichlet
Convergence Lemma}, we find that
\begin{equation*}
\lambda_{i+1}(\Gamma) \leq \|f\|^2_{\Dir}+\varepsilon.
\end{equation*}
We know that $\ip{f}{g}_{\Dir} = \lambda_i(\Gamma) \ip{f}{g}_{L^2}$ for all $g\in
\Zh_{\mu}(\Gamma)$ since $f$ is an eigenfunction of $\Delta$. In particular, setting
$g=f$ tells us that $\|f\|_{\Dir}^2 = \lambda_i(\Gamma)$. Taking $\varepsilon$
sufficiently small provides a contradiction since $\lambda_i(\Gamma) <
\lambda_{i+1}(\Gamma)$. Thus assertion (B) must be true for all $i$.

Now we prove assertion (C). As assertion (B) holds for each fixed $i$, we know that
$\HH_N(i)$ consists of $d=d_i(\Gamma)$ pairwise $\ell^2$-orthogonal functions on $\Gamma$
when $N$ is large. Let $h_N^1, \ldots, h_N^d$ be the elements in $\HH_N(i)$. By normality
of $\HH(i)$, we can pass to a subsequence such that each of these $d$ sequences
$\{h_N^j\}$ converges. Now use Lemma~\ref{Technical Crap 5} to see that the limit
functions must be pairwise $L^2$-orthogonal and have unit $L^2$-norm. The dimension of
the eigenspace of $\Delta$ corresponding to $\lambda_i(\Gamma)$ is $d$, so the
subsequential limits just constructed form a basis for this eigenspace.
\end{proof}

\section{Normality of the family $\HH(i)$} \label{Normality Section}

For the remainder of the paper we assume that assertions (A), (B'), and (C') of the
Reduction Lemma hold up to the case $i-1$, and we prove that they are true for the case
$i$.

\begin{thm}[Eigenvalue Convergence: Part 1] \label{Half The Proof}
We have
\begin{equation*}
\limsup_{N \to \infty} N\lambda_{i,N} \leq \lambda_i(\Gamma).
\end{equation*}
\end{thm}

\begin{proof} Suppose the result is false and pick a subsequence of models
and a positive constant $\delta$ so that $N\lambda_{i,N} > \lambda_i(\Gamma) + \delta$
for all $N$ sufficiently large in our subsequence. Let $f$ be an $L^2$-normalized
eigenfunction for $\Delta$ in $\Zh_{\mu}(\Gamma)$ with eigenvalue $\lambda_i(\Gamma)$. We
now select a further subsequence as in the Approximation Lemma of \S\ref{Reduction Lemma
Section}. For any $\varepsilon
> 0$ and $N$ sufficiently large in our sub-subsequence, we can find $\widetilde{f}_N$ and
apply the Approximation Lemma and Lemma~\ref{Min at Eigenvector} to get
\begin{equation*}
\lambda_i(\Gamma) + \delta < N\lambda_{i,N} = \min_{\substack{g \in \TT_N^{\perp} \\ g
\not\equiv 0}} \frac{\|g\|_{\Dir}^2}{\|g\|_{\ell^2}^2} \leq
\frac{\|\widetilde{f}_N\|_{\Dir}^2}{\|\widetilde{f}_N\|_{\ell^2}^2} < \|f_N\|_{\Dir}^2 +
\varepsilon.
\end{equation*}
We remark that $\TT_N^{\perp}$ is nonempty for $N$ sufficiently large by assertion $(B')$
(which we are assuming holds up to $i-1$). Letting $N$ tend to infinity through our
subsequence and applying Lemma~\ref{Dirichlet Convergence Lemma}, we arrive at the
statement
\begin{equation*} \label{Bound Contradiction}
\lambda_i(\Gamma) + \delta \leq \|f\|_{\Dir}^2 + \varepsilon, \qquad \text{for any
$\varepsilon
> 0$.}\end{equation*}
As $\|f\|^2_{\Dir} = \lambda_i(\Gamma)$, taking $\varepsilon$ sufficiently small produces
a contradiction.
\end{proof}

\begin{cor}  \label{Arzela-Ascoli 1}
The family $\HH(i)$ has uniformly bounded Dirichlet norms.
\end{cor}

\begin{proof}
By definition of an eigenfunction of $Q_N$ with respect to the measure $\mu_N$, we have
immediately that $\|h\|_{\Dir}^2 = N\lambda_{i,N}$ for any $h \in \HH_N(i)$. The previous
theorem shows that $\|h\|_{\Dir}^2$ cannot be arbitrarily large for $h$ lying in
$\HH(i)$.
\end{proof}

\begin{thm}  \label{Arzela-Ascoli 2}
\begin{enumerate}
\item[(a)] The family $\HH(i)$ is equicontinuous on $\Gamma$.
\item[(b)] $\HH(i)$ is uniformly bounded on $\Gamma$.
\item[(c)]$\HH(i)$ is a normal family.
\end{enumerate}
\end{thm}

\begin{proof}
\begin{enumerate}

\item[(a)] Suppose $M$ is a positive real number such that $\|h\|_{\Dir} \leq M$ for all $h \in
\HH(i)$. Now suppose that $\varepsilon >0$ is given and that we have any pair $x,y \in
\Gamma$ with $\dist(x,y) < \varepsilon^2 / M^2$, where $\dist(\cdot, \cdot)$ is the
metric on $\Gamma$. Let $\gamma \subset \Gamma$ be a unit speed path from $x$ to $y$. The
fundamental theorem of calculus and the Cauchy-Schwarz inequality imply that
\begin{equation*}
\begin{aligned}
\left|h(x) - h(y) \right| &= \left| \int_{\gamma} h'(t) dt
\right| \\
&\leq \left( \int_{\gamma} dt \right)^{1/2} \left( \int_{\gamma}
\left|h'(t)\right|^2 dt \right)^{1/2} \\
&\leq \dist(x,y)^{1/2} \|h\|_{\Dir} < \varepsilon.
\end{aligned}
\end{equation*}
The integrals in the above estimate should be viewed as path integrals in which we have
implicitly chosen compatible orientations for consecutive segments along the path
$\gamma$.

\item[(b)] It suffices to show that the family $\HH(i)$ is uniformly bounded on $\bigcup_N V_N$ as this
constitutes a dense subset of $\Gamma$ and all of our functions are continuous.

Recall that if $h \in \HH_N(i)$, then $\int_{\Gamma} h d\mu_N = 0$. By the proof of part
(a) it follows that for any vertex $v \in V_N$ and $h \in \HH_N(i)$,
\begin{equation*}
\begin{aligned}
\left|h(v)\right| &= \left| \int_{\Gamma} \left(h(v) -
h(y)\right) d\mu_N(y)\right| \\
&\leq \int_{\Gamma} \left| h(v) - h(y)\right| d|\mu_N|(y) \\
&\leq \int_{\Gamma} \dist(v,y)^{1/2} \|h\|_{\Dir} d|\mu_N|(y) \\
&\leq \|h\|_{\Dir} |\mu_N|(\Gamma).
\end{aligned}
\end{equation*}
The last inequality follows because the distance between any two points in $\Gamma$ is at
most $1$. Now Corollary~\ref{Arzela-Ascoli 1} and the assumption that the measures
$\mu_N$ have uniformly bounded total variation show that there is a constant $M$,
independent of $N$, so that $|h(v)| \leq M$ for any choice of $h \in \HH_N$ and $v \in
V_N$.

Take any $h \in \HH(i)$ now. Choose $N$ so that $h \in \HH_N(i)$. Note that if $x \in
\Gamma$, then $x$ lies on some segment $e$ of $\Gamma$ with respect to the vertex set
$V_N$. By our choice of $N$, we find $h$ is affine on $e$. Taking $y,z \in V_N$ to be the
endpoints of $e$, we know that there is $t \in [0,1]$ so that $h(x) = (1-t)h(y) + th(z)$.
The triangle inequality implies $|h(x)| \leq \max\{|h(y)|,|h(z)|\} \leq \max_{v_i \in
V_N} |h(v_i)| \leq M$.

\item[(c)] Apply the Arzela-Ascoli theorem using parts (a) and (b) to satisfy the necessary
hypotheses.
\end{enumerate}
\end{proof}

\section{Convergence and approximation results} \label{Approximation Section}

This section contains the proofs of various technical lemmas regarding the integral
operator $\varphi_{\mu}$ and its relation to the discrete integral operator $\p$.

\begin{lem}[\cite{Rud}, p.168, Exercise 16] \label{Rudin Exercise}
Suppose $\{k_N\}$ is an equicontinuous sequence of functions on the metrized graph
$\Gamma$ such that $\{k_N\}$ converges pointwise. Then $\{k_N\}$ converges uniformly on
$\Gamma$.
\end{lem}

\begin{proof}
Fix $\varepsilon > 0$. By equicontinuity there exists $\delta > 0$ such that $|k_N(x) -
k_N(y)| < \varepsilon / 3$ for every $N$ and $x,y \in \Gamma$ with $\dist(x,y) < \delta$.
By compactness, we can cover $\Gamma$ with a finite number of balls of radius $\delta
/2$. Select $y_1, \ldots, y_m \in \Gamma$ so that at least one $y_i$ lies in each ball.

For $j=1, \ldots, m$, let $A_j$ be a positive real number such that $M,N > A_j$ implies
$|k_N(y_j) - k_M(y_j)| < \varepsilon/3$. This is possible because $\{k_N\}$ is pointwise
convergent.  Set $A = \max\{A_j: j=1, \ldots, m\}$.

Suppose $M,N > A$ and that $x \in \Gamma$. Pick $y_i$ so that $\dist(x,y_i) < \delta$.
Then
\begin{equation*}
\left|k_N(x) - k_M(x)\right| \leq \left| k_N(x) - k_N(y_i)\right| + \left| k_N(y_i) -
k_M(y_i)\right| + \left| k_M(y_i) - k_M(x)\right| < \varepsilon.
\end{equation*}
Note that our choice of $M,N$ did not depend on $x$. This shows $\{k_N\}$ is uniformly
Cauchy.
\end{proof}

\begin{lem} \label{Technical Crap 1}
Suppose that $\{k_N\}$ is a sequence of continuous functions on $\Gamma$ which converges
uniformly to a function $k: \Gamma \to \CC$. Then $\varphi_{\mu}(k_N) \to
\varphi_{\mu}(k)$ uniformly with $N$.
\end{lem}

\begin{proof}
Our strategy will be to use Lemma~\ref{Rudin Exercise}. Fix $x \in \Gamma$.  Then, since
$g_{\mu}(x,y)$ is continuous on $\Gamma^2$, it must be uniformly bounded (compactness).
Also, since $\{k_N\}$ converges uniformly to a continuous function $k$, we find that
$k_N(y)$ is uniformly bounded for all $N,y$.  Thus Lebesgue dominated convergence
guarantees that $\varphi_{\mu}(k_N) \to \varphi_{\mu}(k)$ pointwise.

Suppose that $M$ is a positive real number such that $|k_N(y)|\leq M$ for all $N,y$.
Also, since $g_{\mu}$ is continuous on a compact space, it must be uniformly continuous.
Thus, given $\varepsilon > 0$, we can select $\delta >0$ such that $|g_{\mu}(x,y) -
g_{\mu}(x',y)| < \varepsilon/M$ whenever $x,x',y \in \Gamma$ with $\dist(x,x') < \delta$.
Hence for all $x,x',y \in \Gamma$ such that $\dist(x,x') < \delta$, we have
\begin{equation*}
\left| \varphi_{\mu}(k_N)(x) - \varphi_{\mu}(k_N)(x') \right| \leq \int_{\Gamma} \left|
g_{\mu} (x,y) - g_{\mu} (x',y)\right| \cdot \left|k_N(y) \right| dy \leq \varepsilon.
\end{equation*}
This shows the desired equicontinuity.
\end{proof}

\begin{lem} \label{Technical Crap 3.5}
As $N \to \infty$, we have $g_{\mu_N} \to g_{\mu}$ uniformly on $\Gamma^2$.
\end{lem}

\begin{proof}
We intend to apply Lemma~\ref{Rudin Exercise}.

First we show pointwise convergence. By definition of $g_{\mu_N}$ and $g_{\mu}$, it
suffices to show
\begin{itemize} \item $j_{\mu_N}(x,y)
\to j_{\mu}(x,y)$ for each pair $x,y \in \Gamma$, and \item $\int_{\Gamma}
j_{\mu_N}(x,y)d\mu_N(x) \to \int_{\Gamma}j_{\mu} (x,y) d\mu(x)$ for all $y\in \Gamma$.
\end{itemize} The former assertion follows by weak convergence of measures. As for the
second, fix $\varepsilon > 0$ and a basepoint $y_0 \in \Gamma$. As $j_{z}(x,y)$ is a
continuous function on $\Gamma^3$ (a compact space), it is uniformly continuous. So there
exists $\delta > 0$ such that for any $x,x',z \in \Gamma$ with $\dist(x,x')< \delta$, we
have $|j_{z}(x,y_0) - j_{z}(x',y_0)|< \varepsilon$. Now for such $x,x'$, it follows that
\begin{equation*}
\left|j_{\mu_N}(x,y_0) - j_{\mu_N}(x',y_0)\right| \leq \int_{\Gamma} \left|j_{z}(x,y_0) -
j_{z}(x',y_0)\right| d|\mu_N|(z)  \leq \varepsilon |\mu_N|(\Gamma).
\end{equation*}
As the total variations $|\mu_N|(\Gamma)$ are uniformly bounded independent of $N$, we
see $\{j_{\mu_N}(x,y_0)\}$ is equicontinuous in the variable $x$. By Lemma~\ref{Rudin
Exercise} we find $j_{\mu_N}(x,y_0) \to j_{\mu}(x,y_0)$ uniformly in $x$.  For ease of
notation, set $I_{\nu}(\sigma)= \int_{\Gamma} j_{\nu}(x,y_0) d\sigma(x)$ for any pair of
measures $\nu, \sigma$. Now suppose that $N$ is so large that $j_{\mu_N}(x,y_0)$ is
uniformly within $\varepsilon$ of $j_{\mu}(x,y_0)$, and that $I_{\mu}(\mu_N)$ is within
$\varepsilon$ of $I_{\mu}(\mu)$ (by weak convergence). Then
\begin{equation} \label{Pointwise calculation}
\begin{aligned}
\left| I_{\mu_N}(\mu_N) - I_{\mu}(\mu) \right| &\leq \left| I_{\mu_N}(\mu_N) -
I_{\mu}(\mu_N) \right|  + \left| I_{\mu}(\mu_N) - I_{\mu}(\mu)
\right|\\
&\leq \int_{\Gamma} \left|j_{\mu_N}(x,y_0)  - j_{\mu}(x,y_0)\right| \ d|\mu_N|(x) +
\left| I_{\mu}(\mu_N) - I_{\mu}(\mu)
\right|\\
&\leq \varepsilon |\mu_N|(\Gamma) + \varepsilon.
\end{aligned}
\end{equation}
By Proposition~\ref{Doesn't Depend on y}, we see that $\int_{\Gamma} j_{\mu_N} (x,y)
d\mu_N(x)$ and $\int_{\Gamma} j_{\mu}(x,y) d\mu$ are independent of $y$. As $\varepsilon$
was arbitrary, we conclude from \eqref{Pointwise calculation} that $\int_{\Gamma}
j_{\mu_N} (x,y) d\mu_N(x) \to \int_{\Gamma} j_{\mu}(x,y) d\mu$ pointwise for all $y \in
\Gamma$.

All of this allows us to conclude that $g_{\mu_N}(x,y) \to g_{\mu}(x,y)$ pointwise on
$\Gamma^2$. It remains to prove $\{g_{\mu_N}(x,y)\}$ is equicontinuous jointly in $x$ and
$y$.  By uniform continuity of $j$ on $\Gamma^3$, we know there exists $\delta>0$ such
that for all $x,x',y,y',z$ with $\dist\left( (x,y),(x',y')\right) < \delta$ (some metric
on the product space $\Gamma^2$), it must be that $|j_{z}(x,y) - j_{z}(x',y')| <
\varepsilon$. Now for any $x,x',y,y' \in \Gamma$ with $\dist\left( (x,y),(x',y')\right) <
\delta$, we have \begin{align*}\left|g_{\mu_N}(x,y)-g_{\mu_N}(x',y')\right| &=
\left|j_{\mu_N}(x,y)-j_{\mu_N}(x',y')\right|\\  &\leq \int_{\Gamma}
\left|j_{z}(x,y)-j_{z}(x',y')\right| d|\mu_N|(z) \leq \varepsilon|\mu_N|(\Gamma).
\end{align*}
This proves the equicontinuity of $\{g_{\mu_N}(x,y)\}$.
\end{proof}

\begin{lem} \label{Technical Crap 4}
Given any convergent sequence $\{h_N\} \subset \HH(i)$ with limit function $h$ such that
$h_N \in \HH_N(i)$ for each $N$, we find that $\p(h_N) \to \varphi_{\mu}(h)$ uniformly on
$\Gamma$.
\end{lem}

\begin{proof}
We use Lemma~\ref{Rudin Exercise} again. First we show pointwise convergence. Fix $x \in
\Gamma$ and $\varepsilon>0$. As $g_{\mu_N} \to g_{\mu}$ uniformly and $h_N \to h$
uniformly, we can suppose that $N$ is large enough to guarantee $g_{\mu_N}h_N$ is
uniformly within $\varepsilon$ of $g_{\mu}h$ for all $x,y \in \Gamma$. We can also
suppose that $N$ is so large that $\int_{\Gamma} g_{\mu}(x,y)h(y)dy_N$ is within
$\varepsilon$ of $\int_{\Gamma}g_{\mu}(x,y)h(y)dy$. Now we have
\begin{equation*}
\begin{aligned}
\left|\p\{h_N\}(x) - \varphi_{\mu}\{h\}(x)\right| &= \left| \int_{\Gamma}
g_{\mu_N}(x,y)h_N(y)dy_N - \int_{\Gamma}
g_{\mu}(x,y)h(y)dy\right| \\
&\leq \left| \int_{\Gamma} g_{\mu_N}(x,y)h_N(y)dy_N -\int_{\Gamma}
g_{\mu}(x,y)h(y)dy_N \right| \\
& \hspace{0.5in} + \left|\int_{\Gamma} g_{\mu}(x,y)h(y)dy_N -
\int_{\Gamma} g_{\mu}(x,y)h(y)dy \right| \\
&< \int_{\Gamma} \left|g_{\mu_N}(x,y)h_N(y)-g_{\mu}(x,y)h(y)\right| dy_N + \varepsilon
\\ &< 2\varepsilon.
\end{aligned}
\end{equation*}
Hence pointwise convergence holds.

Now we wish to show equicontinuity of the sequence of functions $\{\p(h_N)\}$.  The
sequence $\{g_{\mu_N}(x,y)\}$ is equicontinuous on $\Gamma^2$ by the proof of
Lemma~\ref{Technical Crap 3.5}, and the sequence $\{h_N\}$ is uniformly bounded by some
positive constant $M$. Thus we have
\[\left|\p\{h_N\}(x) - \p\{h_N\}(x')\right| \leq M\int_\Gamma |g_{\mu_N}(x,y) -
g_{\mu_N}(x',y)| dy_N,\] and the integrand can be made arbitrarily small by taking $x$
close to $x'$.
\end{proof}

\begin{lem} \label{Technical Crap 6}
Suppose that $\{h_N\} \subset \HH(i)$ is any convergent sequence with limit function $h$
such that $h_N \in \HH_N(i)$ for each $N$. For each $\varepsilon
> 0$ there exists a real number $N_2$ such that for all $N
> N_2$, we have simultaneously
\begin{itemize}
\item $\|\varphi_{\mu}(h)\|_{L^2} > \|\p(h_N)\|_{\ell^2} - \varepsilon$;
\item $\|h\|_{L^2} < \|h_N\|_{\ell^2} + \varepsilon$.
\end{itemize}
\end{lem}

\begin{proof}
Observe that $h_N \to h$ uniformly, and so $\p(h_N) \to \varphi_{\mu}(h)$ uniformly by
Lemma~\ref{Technical Crap 4}. Now use Lemma~\ref{Technical Crap 5} to assert that for all
$N$ large, $\|h_N\|_{\ell^2}$ is within $\varepsilon$ of $\|h\|_{L^2}$ and
$\|\p(h_N)\|_{\ell^2}$ is within $\varepsilon$ of $\|\varphi_{\mu}(h)\|_{L^2}$.
\end{proof}

\begin{lem}  \label{Technical Crap 3} Suppose $h:\Gamma \to \CC$ is a continuous function such that
$\int_\Gamma h du = 0$ and $h$ is $L^2$-orthogonal to all of the eigenfunctions of
$\Delta$ with associated eigenvalues $\lambda_1(\Gamma), \ldots, \lambda_{i-1}(\Gamma)$.
For each $\varepsilon
> 0$, there exists $\widetilde{h} \in \Zh_{\mu}(\Gamma)$ that satisfies the following two conditions:
\begin{itemize}
\item[(i)] $\widetilde{h}$ is $L^2$-orthogonal to all of the eigenfunctions of $\Delta$ with
associated eigenvalues $\lambda_1(\Gamma), \ldots, \lambda_{i-1}(\Gamma)$; and
\item[(ii)] $\frac{\|\varphi_{\mu}(\widetilde{h})\|_{L^2}}{\|\widetilde{h}\|_{L^2}}
>\frac{\|\varphi_{\mu}(h)\|_{L^2}-\varepsilon}{\|h\|_{L^2}+\varepsilon}$.
\end{itemize}
\end{lem}

\begin{proof}
Using the vertex set $X$, we may decompose $\Gamma$ into a finite number of segments
$\{e_1, \ldots, e_r\}$, each isometric to a closed interval. As $h$ is continuous, we can
apply the Stone-Weierstrass theorem to each segment of $\Gamma$ to get uniform polynomial
approximations of $h|_{e_j}$ with as great an accuracy as we desire. Moreover, we may
define our approximations so that their values agree with the values of $h$ at the
endpoints of each segment, and hence we may glue our approximations together to get a
uniform approximation of $h$. Evidently such approximations must lie in $\Zh(\Gamma)$.

Let us perform this procedure to construct a sequence of functions $\{f_n\} \subset
\Zh(\Gamma)$ that converges uniformly to $h$. Let $\xi_1, \ldots, \xi_p$ be
eigenfunctions of the Laplacian $\Delta$ that form an $L^2$-orthonormal basis for the
direct sum of the eigenspaces associated with the eigenvalues $\lambda_1(\Gamma), \ldots,
\lambda_{i-1}(\Gamma)$. Now set \[H_n(x) = f_n(x) -\sum_{j=1}^p \ip{f_n}{\xi_j}_{L^2}
\xi_j(x)- \int_{\Gamma} f_n d\mu.\] As each $\xi_j$ is orthogonal to constant functions
(constants are eigenfunctions of $\varphi_{\mu}$ with eigenvalue $0$), we see that $H_n
\in \Zh_{\mu}(\Gamma)$ for all $n$ and condition (i) of the lemma is satisfied for all
$H_n$. Since $f_n$ tends uniformly to $h$, we see that the inner products
$\ip{f_n}{\xi_j}_{L^2}$ tend to zero as $n$ tends to infinity. By dominated convergence,
the integral $\int_\Gamma f_n d\mu$ tends to zero as well. Thus $H_n$ converges uniformly
to $h$ on $\Gamma$.

By Lemma~\ref{Technical Crap 1}, we know that $\varphi_{\mu}(H_n) \to \varphi_{\mu}(h)$
uniformly. Since the $L^2$-norm respects this convergence, for any $\varepsilon > 0$ we
may take $n$ sufficiently large and set $\widetilde{h} = H_n$ to obtain condition (ii).
\end{proof}

The eigenvalues of $\varphi_{\mu}$ constitute a sequence of positive numbers tending to
zero. If we label these distinct values by $\alpha_1(\Gamma) > \alpha_2(\Gamma)
>  \alpha_3(\Gamma) > \cdots$, then $\alpha_i(\Gamma) = 1/\lambda_i(\Gamma)$ for
all $i$ (see Theorem~\ref{Reciprocal Continuous Eigenvalues}). We have the following
classical characterization of the $i$th eigenvalue of $\varphi_{\mu}$:

\begin{prop} \label{Technical Crap 7}
Let $\SSS(i)$ denote the span of the eigenfunctions of $\Delta$ with respect to the
measure $\mu$ that are associated to the eigenvalues $\lambda_1(\Gamma), \ldots,
\lambda_{i-1}(\Gamma)$. If $\SSS(i)^{\perp}$ is the $L^2$-orthogonal complement of this
space inside $\Zh_{\mu}(\Gamma)$, then we find that
\begin{equation} \label{Supremum characterization}
\alpha_i(\Gamma) = \sup_{\substack{F \in \SSS(i)^{\perp}  \\
F \not= 0}} \frac{\|\varphi_{\mu}(F)\|_{L^2}}{\|F\|_{L^2}}.
\end{equation}
Moreover, there exists an eigenfunction $F \in \SSS(i)$ which realizes this supremum.
\end{prop}

\begin{proof}
In section~8 of \cite{BR}, it is proved that $\varphi_{\mu}$ is a compact Hermitian
operator on $L^2(\Gamma)$; hence, by the spectral theorem, the eigenfunctions of
$\varphi_{\mu}$ form an orthonormal basis for $L^2(\Gamma)$. In section~15 of \cite{BR}
it is shown that the eigenfunctions of $\varphi_{\mu}$ must lie in $\Zh_{\mu}(\Gamma)$.
We may now identify $\SSS(i)^{\perp}$ with a subspace $W$ of $L^2(\Gamma)$ that is
orthogonal to all of the eigenfunctions of $\Delta$ associated to the eigenvalues
$\lambda_1(\Gamma), \ldots, \lambda_{i-1}(\Gamma)$; moreover, we see that $\varphi_{\mu}$
maps this subspace into itself.

The space $\Zh(\Gamma)$ is dense in $L^2(\Gamma)$ because, for example, smooth functions
are dense in $L^2$. Thus the supremum on the right side of equation~\eqref{Supremum
characterization} is nothing more than the operator norm of $\varphi_{\mu}|_W$. By
\cite{Yng} Theorem~8.10, we find that either $\|\varphi_{\mu}|_W\|$ or
$-\|\varphi_{\mu}|_W\|$ is an eigenvalue of $\varphi_{\mu}$ acting on $W$. The nonzero
eigenvalues of $\varphi_{\mu}$ are all positive, so we may eliminate the latter
possibility. We conclude that there exists an eigenfunction $F \in W$ with associated
eigenvalue $\alpha$ for which $\|\varphi_{\mu}|_W\| = \|\varphi_{\mu}(F)\|_{L^2} /
\|F\|_{L^2} = \alpha$. Evidently $\alpha = \alpha_i(\Gamma)$ since the operator norm is
defined as a supremum and $\alpha_i(\Gamma)$ is the largest eigenvalue of $\varphi_{\mu}$
acting on $\SSS(i)^{\perp}$.
\end{proof}

\section{Proofs of assertions (A), (B'), and (C')} \label{Completion of Proof Section}

Let $\alpha_{1,N} > \alpha_{2,N} > \alpha_{3,N} > \cdots$ denote the nonzero eigenvalues
of the operator $\p$ described in \S\ref{Integral Operator section}. There we saw that
for each $i$, the eigenvalues of $\p$ are related to the eigenvalues of $Q_N$ with
respect to the measure $\mu_N$ via $N \lambda_{i,N} = 1/\alpha_{i,N}$.

Recall that we are assuming assertions (A), (B'), and (C') of the Reduction Lemma hold up
to the case $i-1$.

\begin{thm}[Eigenvalue Convergence: Part 2] \label{Eigenvalue Convergence Part 2}
We have \begin{equation*} \lambda_i(\Gamma) \leq \liminf_{N \to \infty} N\lambda_{i,N}.
\end{equation*}
\end{thm}

\begin{proof}
Suppose the theorem is false. As $N\lambda_{i,N} = 1/\alpha_{i,N}$ and $\lambda_i(\Gamma)
= 1/\alpha_i(\Gamma)$, our supposition is equivalent to $\alpha_i(\Gamma) < \limsup_{N
\to \infty} \alpha_{i,N}$. We may pick a subsequence of $\{G_N\}$ and a positive constant
$\eta$ so that $\alpha_i(\Gamma) + \eta < \alpha_{i,N}$ for each $N$ in the subsequence.
For each such $N$, we also pick $h_N \in \HH_N(i)$. By normality of $\HH(i)$, upon
passage to a further subsequence we may assume that $\{h_N\}$ tends uniformly to some
continuous function $h$ as $N$ goes to infinity through this subsequence.

Fix $\varepsilon > 0$ and select $\widetilde{h}$ as in Lemma~\ref{Technical Crap 3}. For
$N$ large, the conclusions of Lemma~\ref{Technical Crap 6} are true. It now follows that
for any $N$ sufficiently huge in our subsequence, we have
\begin{equation*}
\begin{aligned}
\alpha_i(\Gamma) &= \sup_{\substack{F \in \SSS(i)^{\perp}  \\
F \not= 0}} \frac{\|\varphi_{\mu}(F)\|_{L^2}}{\|F\|_{L^2}}, \quad \textrm{by
Lemma~\ref{Technical Crap 7}} \\
&\geq \frac{\|\varphi_{\mu}(\widetilde{h})\|_{L^2}} {\|\widetilde{h}\|_{L^2}}, \quad
\textrm{since
$\widetilde{h} \in \SSS(i)^{\perp}$} \\
&> \frac{\|\varphi_{\mu}(h)\|_{L^2} -\varepsilon}{\|h\|_{L^2}+\varepsilon}, \quad
\textrm{by
Lemma~\ref{Technical Crap 3}} \\
&\geq \frac{\|\p(h_N)\|_{\ell^2} -2\varepsilon}{\|h_N\|_{\ell^2}+2\varepsilon}, \quad
\textrm{by
Lemma~\ref{Technical Crap 6}} \\
&=\frac{\alpha_{i,N}\|h_N\|_{\ell^2}
-2\varepsilon}{\|h_N\|_{\ell^2}+2\varepsilon}, \quad
\textrm{since $h_N$ is an eigenfunction of $\p$} \\
&= \frac{\alpha_{i,N} - 2\varepsilon}{1+2\varepsilon}\\
&> \frac{\alpha_i(\Gamma) + \eta - 2\varepsilon}{1+2\varepsilon}.
\end{aligned}
\end{equation*}
This furnishes us with a contradiction when $\varepsilon$ is sufficiently small.
\end{proof}

Clearly the amalgamation of Theorems~\ref{Half The Proof}~and~\ref{Eigenvalue Convergence
Part 2} prove that assertion (A) of our Main Theorem holds for $i$. As for assertion (B')
and the remainder of assertion (C') in the Reduction Lemma, we have the following two
corollaries.

\begin{cor}\label{Limits are Eigenfunctions}
The subsequential limits of $\HH(i)$ are eigenfunctions of $\Delta$ with associated
eigenvalue $\lambda_i(\Gamma)$. The limits have unit $L^2$-norm and lie in
$\Zh_{\mu}(\Gamma)$.
\end{cor}

\begin{proof}
Let $\{h_N\}$ be a convergent subsequence of $\HH(i)$ with limit function $h$. We know
that $\p(h_N)(x) = \frac{1}{N\lambda_{i,N}} h_N(x)$ for all $x  \in \Gamma$. Now let $N$
go to infinity and apply Lemma~\ref{Technical Crap 4} and assertion (A) to see that
$\varphi_{\mu}(h) = \frac{1}{\lambda_i(\Gamma)} h$. Hence $h$ is an eigenfunction of
$\varphi_{\mu}$. By Theorem~12.1 and Proposition~15.1 in \cite{BR} we see that $h$ is an
eigenfunction of $\Delta$ lying in $\Zh_{\mu}(\Gamma)$. Lemma~\ref{Technical Crap 5}
shows $\|h\|_{L^2}=1$.
\end{proof}

\begin{cor} \label{Eigenspace Stabilization}
For $N$ sufficiently large, it is true that $d_{i,N} \leq d_i(\Gamma)$.
\end{cor}

\begin{proof} Set $d=d_i(\Gamma)$. Suppose that that $d_{i,N} > d$ for some infinite
subsequence of $\{G_N\}$. For each such $N$, select distinct functions $h_N^1, \ldots,
h_N^{d+1} \in \HH_N(i)$. Using the normality of $\HH(i)$, we may pass to a further
subsequence and assume that there are $d+1$ limit functions $h^1, \ldots, h^{d+1}$, with
$h_N^j \to h^j$ uniformly along this subsequence for each $j$. By Corollary~\ref{Limits
are Eigenfunctions}, we find that each $h^j$ is an eigenfunction of $\Delta$ with
associated eigenvalue $\lambda_i(\Gamma)$. But the functions $h_N^1, \ldots, h_N^{d+1}$
are $\ell^2$-orthogonal for each $N$, and Lemma~\ref{Technical Crap 5} implies that the
limit functions must be $L^2$-orthogonal. As the eigenspace corresponding to the
eigenvalue $\lambda_i(\Gamma)$ has dimension $d$, this is an absurdity.
\end{proof}

We have now exhibited the induction step in our proof of the Main Theorem; i.e., the
proof is complete.

\section{Acknowledgements}

This work could not have been completed without the gracious support of the NSF VIGRE
grant at the University of Georgia or the support of the Columbia University mathematics
department. I would like to thank Matthew Baker and Robert Rumely for introducing me to
the subject and for their most helpful guidance toward completing this project. I would
also like to thank David Glickenstein for pointing me toward the relevant work being done
in China and Japan. The members of the 2003 UGA REU on metrized graphs deserve a round of
applause for their hard work and enthusiasm. My wife has always provided me with a loving
environment in which to work, and for this she should be praised most highly.


\vspace{0.3 in}

\begin{thebibliography}{[1]}


    \bibitem[BF]{BF}{Matt Baker and Xander Faber. Metrized graphs,
    electrical networks, and Fourier analysis. Submitted. }

    \bibitem[Bo]{Bo}{B{\'e}la Bollob{\'a}s, \textit{Modern Graph
    Theory}, Springer-Verlag, New York, 1998.}

    \bibitem[BR]{BR}{Matthew Baker and Robert Rumely. Harmonic analysis on
    metrized graphs. To appear in \textit{Canad. J. Math.}}


    \bibitem[CR]{CR}{Ted Chinburg and Robert Rumely. The capacity
    pairing. \textit{J. reine angew. Math.} \textbf{434} (1993), 1-44.}

    \bibitem[KoFu1]{KoFu1}{Koji Fujiwara. Convergence of the eigenvalues of Laplacians in
    a class of finite graphs. \textit{Geometry of the spectrum}. (Seattle, WA, 1993)
    Contemp. Math., \textbf{173}, Amer. Math. Soc., Providence, RI, 1994. 115-120.}

    \bibitem[KoFu2]{KoFu2}{Koji Fujiwara. Eigenvalues of Laplacians on a closed Riemannian manifold and its
    nets. \textit{Proc. Amer. Math. Soc.} \textbf{123} (1995), no. 8, 2585-2594.}

    \bibitem[KeFu]{KeFu}{Kenji Fukaya. Collapsing of Riemannian manifolds and eigenvalues of Laplace
    operator. \textit{Invent. Math.} \textbf{87} (1987), no. 3, 517-547.}



    \bibitem[Ku]{Ku}{Peter Kuchment. Quantum graphs: I. Some basic
    structures. \textit{Waves Random Media.} \textbf{14} (2004), S107-S128.}

    \bibitem[Mo]{Mo}{Bojan Mohar. The Laplacian Spectrum of
    Graphs.  \textit{Graph Theory, Combinatorics, and
    Applications.} Vol. 2, Wiley, 1991, 871-898.}

    \bibitem[Rud]{Rud}{Walter Rudin. \textit{Principles of Mathematical
    Analysis}. 3rd ed., McGraw-Hill, 1976, New York.}

    \bibitem[Yng]{Yng}{Nicholas Young. \textit{An Introduction to Hilbert Space}.
    Cambridge University Press, 1988, New York.}

\end{thebibliography}
\end{document}